\documentclass[reqno,12pt]{amsart}
\usepackage{amssymb}
\usepackage{amsmath}

\usepackage[all]{xy}

\def\id{\mathop{\rm id}\nolimits}

\def\Ad{\mathop{\rm Ad}\nolimits}

\def\ev{\mathop{\rm ev}\nolimits}

\def\be{\begin{equation}}
\def\ee{\end{equation}}

\newtheorem{thm}{Theorem}[section]
\newtheorem{lem}[thm]{Lemma}

\newtheorem{cor}[thm]{Corollary}

\newtheorem{rmk}[thm]{Remark}

\def\co{{\text{co}}}

\date{}
\author{V. Manuilov and K. Thomsen}
\title{$E$-theory is a special case of $KK$-theory}

\begin{document}

\maketitle

 \begin{abstract}
Let $A$ and $B$ be $C^*$-algebras, $A$ separable, and $B$
$\sigma$-unital and stable. It is shown that there are natural
isomorphisms
\begin{center}
 $
E(A,B) = KK(SA,Q(B)) = [SA,Q(B) \otimes \mathbb K],
 $
\end{center}

\noindent
where $SA=C_0(0,1)\otimes A$, $[\cdot,\cdot]$ denotes the set of homotopy
classes of $*$-homomorphisms, $Q(B) = M(B)/B$ is the generalized Calkin
algebra and $\mathbb K$ denotes the $C^*$-algebra of compact operators of
an infinite dimensional separable Hilbert space.

\end{abstract}

\section{Introduction}\label{s1}

Connes and Higson have introduced a variant of Kasparov's $KK$-theory,
\cite{CH}, \cite{K}, called $E$-theory, originally with the purpose of
realizing
a theory developed by Higson, \cite{H2}, which is half-exact with respect
to arbitrary extensions of (separable) $C^*$-algebras. $E$-theory relates
directly to Kasparov's theory via a natural map
$$
KK(A,B) \to E(A,B),
$$
which is an isomorphism when $A$ is nuclear, but which is not always
injective, reflecting the fact that $KK$-theory is not half-exact in
general, \cite{S}. By specializing, this map gives us also a natural map
$$
KK(SA,Q(B \otimes \mathbb K)) \to E(SA, Q(B \otimes \mathbb K)),
$$
where $SA = C_0(0,1) \otimes A$ is the suspension of $A$, and $Q(B \otimes
\mathbb K) = M(B\otimes \mathbb K)/B\otimes \mathbb K$ is the generalized
Calkin algebra, or `corona' algebra, of the stabilized $C^*$-algebra
$B\otimes \mathbb K$. Thanks to the unrestricted excision properties of
$E$-theory, combined with Bott-periodicity and stability, there is a
natural isomorphism $E(SA,Q(B \otimes \mathbb K)) = E(A,B)$, and we have
therefore a natural map
\begin{equation}\label{1!}
KK(SA,Q(B \otimes \mathbb K)) \to E(A, B) .
\end{equation}
The main result of the present paper is that this map is an isomorphism,
when $A$ is separable and $B$ $\sigma$-unital. The method of proof is easily
explained. In \cite{MT3} we introduced a method which produces a genuine
$*$-homomorphism $\widetilde{\varphi}: A \to Q(B \otimes \mathbb K)$ out of an
asymptotic homomorphism $\varphi = \left( \varphi_t\right)_{t \in [1,\infty)}
: A \to Q(B)$. We will show that this construction gives rise to a group
homomorphism
$$
[[SA, Q(B) \otimes \mathbb K]] \to [SA, Q(B \otimes \mathbb K)
\otimes \mathbb K],
$$
where brackets $[\cdot,\cdot]$ (resp. double brackets $[[\cdot,\cdot]]$)
denote the set of homotopy classes of homomorphisms (resp. of asymptotic
homomorphisms).
Composed with the obvious map $[SA, Q(B \otimes \mathbb K) \otimes \mathbb
K] \to KK(SA, Q(B \otimes \mathbb K))$, we obtain a natural map $E(A,B) =
 [[SA, Q(B) \otimes \mathbb K]] \to KK(SA, Q(B \otimes \mathbb K))$, which
we show is an inverse for the map (\ref{1!}). As a by-product of the proof
we also obtain the description of $E$-theory as homotopy classes of
$*$-homomorphisms spelled out in the abstract.

Despite the simplicity in the idea of proof, the actual realization of the
approach is quite technical, and it will occupy the remaining part of the
paper.

It should be pointed out that an alternative translation between
$KK$-theory and $E$-theory was obtained by the second-named author in
\cite{T}; specifically, it was shown that
$$
E(A,B) = KK(A, C_b([1,\infty), B \otimes \mathbb K)/C_0([1,\infty),
B \otimes \mathbb K)),
$$
when $A$ and $B$ are both separable. However, $ C_b([1,\infty), B
\otimes \mathbb K)/C_0([1,\infty),
B \otimes \mathbb K)$ is a nasty $C^*$-algebra which is never
$\sigma$-unital (unless $B = 0$), and is certainly a rather unnatural
gadget to consider in relation to $KK$-theory. In contrast $Q(B\otimes
\mathbb K)$ is always unital and a 'must' in every theory dealing with
extensions of $C^*$-algebras.

In \cite{MT3} we came very close to the
conclusion that $E(A,B) = [SA, Q(B) \otimes \mathbb K]$, in that we proved
the equality $E(A,B) = \varinjlim_n [SA, Q(B \otimes \mathbb K) \otimes M_n]$.
By combining this with the result in the present paper we see that
$\varinjlim_n [SA, Q(B \otimes \mathbb K) \otimes M_n] = [SA, Q(B \otimes
\mathbb K) \otimes \mathbb K]$. This equality seems so plausible that one
is tempted to write down a direct proof. However, we are unaware of any
proof which avoids the use of asymptotic homomorphisms.

\section{Fundamental notation and terminology}\label{s2}

Let $A$ and $B$ be $C^*$-algebras, $M(B)$ the multiplier algebra of $B$ and
$q_B : M(B) \to M(B)/B = Q(B)$ the quotient map. An asymptotic homomorphism
$\varphi = (\varphi_t)_{t\in[1,\infty)}:A\to Q(B)$ will also be called
\emph{an asymptotic extension} of $A$ by $B$. Two asymptotic extensions
$\varphi, \psi :A\to Q(B)$ are \emph{unitarily equivalent} when there is a
normcontinuous path
$(u_t)_{t \in [1, \infty)}$, of unitaries in $M(B)$ such
that $\lim_{t \to \infty} \Ad q_B(u_t) \circ \psi_t(a) - \varphi_t(a) = 0$
for all $a \in A$. An asymptotic extension $\varphi: A \to Q(B)$ is
said to be {\it split} (as an asymptotic extension) when there exists an
asymptotic homomorphism $(\lambda_t)_{t\in[1,\infty)}:A\to M(B)$ such that
$\lim_{t \to \infty} q_B\circ \lambda_t(a) -\varphi_t(a) = 0$ for all $a
\in A$. (We will show along the way that this is equivalent to the apparently
stronger condition that there exists an asymptotic homomorphism
$(\lambda'_t)_{t\in[1,\infty)}:A\to M(B)$ such that $q_B\circ \lambda'_t(a)=
\varphi_t(a)$ for all $a\in A$, $t\in[1,\infty)$,
provided only that $\varphi$ is equi-continuous,
cf. Lemma \ref{crux2}.) An asymptotic homomorphism $\varphi = (\varphi_t
)_{t\in[1,\infty)}: A \to B$ is \emph{equi-continuous} when the family of
maps $(\varphi_t)_{t \in [1,\infty)} : A \to B$ is an equi-continuous
family,
i.e. has the property that for all $a \in A$ and all $\epsilon > 0$ there
is a $\delta > 0$ such that $\|\varphi_t(a) - \varphi_t(b)\| \leq \epsilon$
for all $t \in [1,\infty)$ and all $b \in A$ with $\|a-b\| \leq \delta$. It
is well-known that any asymptotic homomorphism $\varphi : A \to B$ is
asymptotically identical with an equi-continuous asymptotic homomorphism
$\psi$, in the sense that $\lim_{t \to \infty} \varphi_t(a) - \psi_t(a) = 0$
for all $a \in A$. Thanks to this it is almost always possible to restrict
the attention to asymptotic homomorphisms that are equi-continuous. This is
very useful, also for asymptotic extensions, and to make it even more useful
we make the following observation.

\begin{lem}\label{obs} Let $\varphi : A \to Q(B)$ be an equi-continuous
asymptotic extension. It follows that there is an equi-continuous family of
maps $(\widehat{\varphi}_t)_{t \in [1,\infty)} :
A \to M(B)$ such that $q_B \circ
\widehat{\varphi}_t = \varphi_t$ for all $t \in [1,\infty)$ and $\sup_{t \in [
1,\infty)} \|\widehat{\varphi}_t(a)\| < \infty$ for all $a \in A$.
\end{lem}
\begin{proof} Since $\varphi$ is an asymptotic homomorphism,
$\limsup_{t \to \infty} \|\varphi_t(a)\| \leq \|a\|$ for all $a \in A$.
Since $\varphi$
is equi-continuous, we can define a continuous map $\Phi : A \to
C_b([1,\infty), Q(B))$ such that $\Phi(a)(t) = \varphi_t(a)$. There is a
surjective $*$-homomorphism $Q : C_b([1,\infty), M(B)) \to
C_b([1,\infty), Q(B))$ such
that $Q(f)(t) = q_B(f(t))$. By the Bartle-Graves selection theorem there
is a continuous section $S : C_b([1,\infty), Q(B)) \to  C_b([1,\infty),
M(B))$ for $Q$, and we set $\widehat{\varphi}_t(a) = \left(S \circ
\Phi(a)\right)(t)$.
\end{proof}

We will refer to $\widehat{\varphi} = \left( \widehat{\varphi}_t\right)_{t \in [1
,\infty)}$ as \emph{an equi-continuous lift} of $\varphi$. Recall, cf.
\cite{MT1}, that a \emph{discretization},
$\{\varphi_{t_n}\}_{n =1}^{\infty}$, of an asymptotic
homomorphism $\varphi = \left( \varphi_t\right)_{t \in [1,\infty)} :
A \to B$ is given by a sequence $t_1 < t_2 < t_3 < \dots$ in $[1,\infty)$
such that
\begin{enumerate}
\item[(d1)] $\lim_{n \to \infty} t_n = \infty$, and
\item[(d2)] $\lim_{n \to \infty} \sup_{t \in [t_n, t_{n+1}]} \|\varphi_t(a)
- \varphi_{t_n}(a)\| = 0$ for all $a \in A$.
\end{enumerate}
When $A$ is separable discretizations always exist, cf. \cite{MT1}.

For any $C^*$-algebra $A$ we denote by $SA$ the suspension of $A$, i.e. $SA=
C_0(0,1) \otimes A$, and we denote by $IA$ the $C^*$-algebra $C[0,1]
\otimes A$. For any $s \in [0,1]$, $\ev_s : IA \to A$ is then the
$*$-homomorphism given by evaluation at $s$, i.e. $\ev_s(f) = f(s)$.

Concerning extensions we will adopt the terminology from \cite{MT3}; in
particular, we will call an extension $\varphi : A \to Q(B)$
\emph{asymptotically split} when there is an asymptotic homomorphism
$\pi : A \to M(B)$ such
that $\varphi(a) = q_B \circ \pi_t(a)$ for all
$a\in A$, $t\in[1,\infty)$. Concerning asymptotic
homomorphism we shall use the standard notation and terminology; in
particular $[[A,B]]$ will denote the homotopy classes of asymptotic
homomorphisms from $A$ to $B$, and $[[A,B]]_{cp}$ will denote the homotopy
classes of
completely positive asymptotic homomorphsims from $A$ to $B$, where we
call an asymptotic homomorphism $\varphi = \left( \varphi_t\right)_{t \in
[1,\infty)} : A \to B$ \emph{completely positive} when the individual maps
$\varphi_t$ are all completely positive contractions. It is fundamental to
our approach that KK-theory can be realized by using completely positive
asymptotic homomorphisms, specifically that $KK(A,B) = [[SA, SB \otimes
\mathbb K]]_{cp}$, cf. \cite{H-LT}.



\section{The basic construction}\label{s3}

Let $D$ and $E$ be $C^*$-algebras, $D$ separable, $E$ $\sigma$-unital. Let
$(\varphi_t)_{t\in[1,\infty)}:D\to Q(E)$ be an equi-continuous asymptotic
extension. We will construct an extension $\varphi^f$ of $D$ by $E$ out
of $\varphi$. The construction was introduced and used in \cite{MT3}, and
uses Voiculescu's tri-diagonal projection trick from \cite{V}.
Let $b$ be a strictly positive element in $E$ of norm $\leq 1$. A {\it unit
sequence} (cf. \cite{MT3})
in $E$ is a sequence $\{u_n\}_{n=0}^{\infty} \subseteq
E$ such that
\begin{enumerate}
\item[(u1)] for each $n=0,1,2,\dots$
there is a continuous function $f_n : [0,1] \to [0,1]$ which
is zero in a neighbourhood of $0$ and $u_n = f_n(b)$,
\item[(u2)] $u_{n+1}u_n = u_n $ for all $n$, and
\item[(u3)] $\lim_{n \to \infty} u_n x = x$ for any $x \in E$.
\end{enumerate}
Unit sequences exist by elementary spectral theory. Given a unit sequence
$\{u_n\}_{n=0}^{\infty}$ we set $\Delta_0 = \sqrt{u_0}$ and $\Delta_j =
\sqrt{u_j - u_{j-1}}$, $j \geq 1$. Note that (u2) implies that
\begin{equation}\label{BB}
\Delta_i\Delta_j = 0, \ |i-j|\geq 2.
\end{equation}

\begin{lem}\label{E1} For any norm-bounded sequence $\{m_j\} \subseteq M(E)$,
and any $k \in \mathbb N$, the sum $\sum_{j=0}^{\infty} \Delta_j m_j
\Delta_{j+k}$ converges in the strict topology to an element of $M(E)$,
$$
\|\sum_{j=0}^{\infty} \Delta_j m_j \Delta_{j+k} \| \leq \sup_j \|m_j\|
$$
and
$$
\|q_E\Bigl(\sum_{j=0}^{\infty} \Delta_j m_j \Delta_{j+k}\Bigr)\| \leq
\limsup_n \|m_n\|.
$$
In particular, $\sum_{j=0}^{\infty} \Delta_j m_j \Delta_{j+k} \in E$ when
$\lim_{j \to \infty} \|m_j\| = 0$.
\end{lem}
\begin{proof} For any $e \in E$ and any $N < M$ in $\mathbb N$ we have the
estimates
$$
\|\sum_{j=N}^M \Delta_jm_j\Delta_{j+k}e\|^2 \leq \|\sum_{j=N}^M
\Delta_jm_jm_j^*\Delta_{j}\|\|\sum_{j=N}^M e^*\Delta_{j+k}^2e\| \leq
\|\sum_{j=N}^M e^*\Delta_{j+k}^2e\| \max_{N \leq i\leq M}\|m_i\|^2,
$$
and
$$
\|\sum_{j=N}^M \Delta_{j+k}m_j^*\Delta_{j}e\|^2 \leq \|\sum_{j=N}^M
\Delta_{j+k}m_j^*m_j\Delta_{j+k}\|\|\sum_{j=N}^M e^*\Delta_{j}^2e\| \leq
\|\sum_{j=N}^M e^*\Delta_{j}^2e\| \max_{N \leq i\leq M}\|m_i\|^2 .
$$
It follows that $\sum_{j=0}^{\infty} \Delta_j m_j \Delta_{j+k}$ converges
in the strict topology, as asserted. We can then take $N = 0$ and let $M$
tend to infinity in the first estimate to conclude that
$$
\|\sum_{j=0}^{\infty} \Delta_j m_j \Delta_{j+k} \| \leq \sup_j \|m_j\| .
$$
If we only let $M$ go to infinity we see that
$$
\|q_E\Bigl(\sum_{j=0}^{\infty} \Delta_j m_j \Delta_{j+k}\Bigr)\| \leq
\|\sum_{j=N}^{\infty} \Delta_j m_j \Delta_{j+k}\| \leq \sup_{i \geq N}
\|m_i\|
$$
for all $N$.
\end{proof}

Let $(\widehat{\varphi}_t)_{t\in[1,\infty)} : D \to M(E)$
be an equi-continuous lift of $\varphi$.
We will say that the pair
$$\left( \left(\widehat{\varphi}_t\right)_{t \in [1,\infty)},
\{u_n\}_{n = 0}^{\infty}\right)
$$
is a \emph{compatible pair} for $\varphi$ when
\begin{enumerate}
\item[(p1)]  $\lim_{n \to \infty} \sup_{t \in [1, n+3]}
\|u_n\widehat{\varphi}_t (a) - \widehat{\varphi}_{t}(a)u_n\| = 0$
\end{enumerate}
for all $a \in D$.

\begin{lem}\label{ex} For any unit sequence $\{v_n\}_{n = 0}^{\infty}$ in
$E$ there is a unit sequence $\{u_n\}_{n =0}^{\infty}$ such that $u_n$
is in the convex hull $\co \{v_i : i \geq n\}$ for all $n$, and
$$\left( \left(\widehat{\varphi}_t\right)_{t \in [1,\infty)},
\{u_{n}\}_{n=0}^{\infty}\right)
$$
is a compatible pair.
\end{lem}
\begin{proof} Let $F_1 \subseteq F_2 \subseteq F_3 \subseteq \dots$ be a
sequence of finite sets with dense union in $D$.
Since $(\widehat{\varphi}_t)_{t\in [1,\infty)}$
is an equi-continuous family of maps it suffices to
construct a sequence $\{u_n\}_{n =0}^{\infty}$ such that
\begin{enumerate}
\item[(v1)] $\|u_{n} \widehat{\varphi}_{t}(a) -
\widehat{\varphi}_{t}(a)u_{n}\| \leq
\frac{1}{n}$ for all $ a \in F_n$ and all $t \in [1,n+1]$,
\item[(v2)] $u_{n}u_{n -1} = u_{n-1}, n \geq 1$,
\item[(v3)] $u_n \in \co\{v_i : i \geq n\}$,
\end{enumerate}
for all $n$. The $u_n$'s are constructed inductively. Since the argument for
the induction start is contained in the
argument for the induction step, we will only give the latter. So assume
that $u_{n-1}$ has been found. Since $u_{n-1} \in \co \{v_i\}$, there is an
$N \geq n$ such that $v_i u_{n-1} = u_{n-1}$ for all $i \geq N$. Since
$\co \{v_i : i \geq N\}$ is a convex approximate unit in $E$, the existence
of $u_{n}$ follows from \cite{A} or \cite{P}.
\end{proof}

\begin{lem}\label{ex0} Let $\left(\widehat{\varphi}_t\right)_{t \in [1,\infty)}$
be an equi-continuous lift of $\varphi$. There exists a sequence $t_1 < t_2
< t_3 < \dots$ in $[1,\infty)$ such that $\{\varphi_{t_n}\}_{n=1}^{\infty}$
is a discretization of $\varphi$ and
\begin{enumerate}
\item[(t1)] $\lim_{n \to \infty} \sup_{t \in [t_n,t_{n+1}]}
\|\widehat{\varphi}_t(a) - \widehat{\varphi}_{t_n}(a)\| = 0$ for all $a \in D$, and
\item[(t2)] $t_n \leq n$ for all $n \in \mathbb N$.
\end{enumerate}
\end{lem}
\begin{proof} Let $F_1 \subseteq F_2 \subseteq F_3 \subseteq \dots$ be a
sequence of finite sets with dense union in $D$. By uniform continuity there
are $N_n$'s in $\mathbb N$, $N_n \geq 1$, such that,
for any $s,t \in [n,n+1]$,
$$
|s -t| \leq \frac{1}{N_n}  \Rightarrow \|\widehat{\varphi}_t(a)-
\widehat{\varphi}_s(a)\| \leq \frac{1}{n}
$$
for all $a \in F_n$. Set $N_0 = 0$
and $t_j = k + \frac{j -  (N_0 + N_1 + \dots + N_k)}{N_{k+1}}$ when
$N_0 + N_1 + \dots + N_k \leq j \leq  N_0 + N_1 + \dots + N_k + N_{k+1}$.
Then $\lim_{n \to \infty} \sup_{t \in [t_n, t_{n+1}]} \|\widehat{\varphi}_t(a) -
\widehat{\varphi}_{t_n}(a)\| = 0$ for all $a \in \bigcup_n F_n$. By
equi-continuity of $(\widehat{\varphi}_t)_{t \in [1,\infty)}$ the same
must hold for all $a \in D$. Since $\lim_{n \to \infty} t_n = \infty$,
$\{\varphi_{t_n}\}_{n =1}^{\infty}$ is a discretization of $\varphi$.
Note that $t_n \leq n$ for all $n$ by construction.
\end{proof}

\begin{lem}\label{E5}
Let $\left( \left(\widehat{\varphi}_t\right)_{t \in [1,\infty)},
\{u_n\}_{n = 0}^{\infty}\right)$ be a compatible pair for $\varphi$.
There is then a sequence $n_0 < n_1 < n_2 < \dots $ in $\mathbb N$ such that
\begin{enumerate}
\item[(t3)] $n_i - n_{i-1} > i + 1$ for all $i \geq 1$,
\item[(t4)] $\lim_{i \to \infty} \sup_{j \geq n_i} \sup_{t \in [1,i+3]}
\|\left(1-u_j\right)\left( \widehat{\varphi}_t(a)\widehat{\varphi}_t(b) -
\widehat{\varphi}_t(ab)\right)\| - \|\varphi_t(a)\varphi_t(b) - \varphi_t(ab)\| =
0$,
\item[(t5)] $ \lim_{i \to \infty} \sup_{j \geq n_i} \sup_{t \in [1, i+3]} \|
\left(1 -u_j\right) \left( \widehat{\varphi}_t(a) + \lambda \widehat{\varphi}_t(b)
- \widehat{\varphi}_t(a + \lambda b)\right)\| -  \|\varphi_t(a) + \lambda
\varphi_t(b) - \varphi_t(a + \lambda b)\| = 0$,
\item[(t6)]  $ \lim_{i \to \infty} \sup_{j \geq n_i} \sup_{t \in [1, i+3]}
\|\left(1-u_j\right)\left(\widehat{\varphi}_t(a^*) - \widehat{\varphi}_t(a)^*\right)
\| -  \|\varphi_t(a^*) - \varphi_t(a)^*\| = 0$
\end{enumerate}
for all $a,b \in D$ and all $\lambda \in \mathbb C$.
\end{lem}

\begin{proof} Let $F_1 \subseteq F_2 \subseteq F_3 \subseteq \dots$ be a
sequence of finite sets with dense union in $D$, and note that
\begin{eqnarray*}
K_i&=&
\{ \widehat{\varphi}_t(a)\widehat{\varphi}_t(b) -  \widehat{\varphi}_t(ab) :
t \in [1,i+3], \ a,b \in F_i\},\\
&\cup&\{ \widehat{\varphi}_t(a) + \lambda
\widehat{\varphi}_t(b) - \widehat{\varphi}_t(a + \lambda b) : t \in [1,i+3],
\ a,b \in F_i, \ \lambda \in \mathbb C,\ |\lambda| \leq i\}\\
&\cup&\{\widehat{\varphi}_t(a^*) - \widehat{\varphi}_t(a)^* :
t \in [1,i+3],\ a \in F_i\}
\end{eqnarray*}
is a compact subset of $M(E)$. Since $\{u_n\}_{n=0}^{\infty}$
is an approximate unit for $E$ there is an $l_i \in
\mathbb N$ so large that
$$
\|\left(1 - u_j\right) m\| \leq \|q_E(m)\| + \frac{1}{i+1}
$$
for all $m \in K_i$ and all $j \geq l_i$. Set $n_0 = l_0$ and $n_i =
\max \{l_i, n_{i-1} + i +2\}$, $i \geq 1$. Then (t3) holds, and
(t4)-(t6) hold for
all $a,b \in \bigcup_n F_n$ and all $\lambda \in \mathbb C$. We can
therefore conclude by appealing to the equi-continuity of
$(\widehat{\varphi}_t)_{t \in [1,\infty)}$ and $(\varphi_t)_{t \in
[1,\infty)}$.
\end{proof}

Let $\left( \left(\widehat{\varphi}_t\right)_{t \in [1,\infty)},
\{u_n\}_{n = 0}^{\infty}\right)$ be a compatible pair for $\varphi$.
Let $\{\varphi_{t_n}\}_{n = 1}^{\infty}$ be a discretization of
$(\varphi_t)_{t \in [1,\infty)}$ such that (t1) and (t2) hold. Let
$\{n_i\}_{i=0}^{\infty}$ be any sequence in $\mathbb N$ such that (t3)-(t6)
hold. In this setting, put
$$
\Delta_0 = \sqrt{u_{n_0}}
$$
and
$$
\Delta_j = \sqrt{u_{n_j} - u_{n_{j-1}}}, \ \ \ \ \ \ \ j \geq 1.
$$

\begin{lem} There is then a $*$-homomorphism $\varphi^f : D \to Q(E)$ such
that
$$
\varphi^f(d) = q_E\Bigl(\sum_{j=0}^{\infty} \Delta_j
\widehat{\varphi}_{t_{j+1}}(d) \Delta_j\Bigr)
$$
for all $d \in D$.
\end{lem}
\begin{proof} Since $\{\|\widehat{\varphi}_t(d)\| : t \in [1, \infty)\}$ is
bounded, we can apply Lemma \ref{E1} to see that $\sum_{j=0}^{\infty}
\Delta_j \widehat{\varphi}_{t_{j+1}}(d) \Delta_j$ converges in the strict topology.
To check that $\varphi^f$ is a $*$-homomorphism, we calculate modulo $E$,
with the aid of Lemma \ref{E1}:
\begin{eqnarray*}
&& \Bigl(\sum_{j=0}^{\infty} \Delta_j \widehat{\varphi}_{t_{j+1}}(a)
\Delta_j \Bigr) \Bigl(\sum_{j=0}^{\infty} \Delta_j
\widehat{\varphi}_{t_{j+1}}(b)\Delta_j\Bigr) \\
&& = \sum_{j=0}^{\infty} \Delta_j \widehat{\varphi}_{t_{j+1}}(a) \Delta_j^2
\widehat{\varphi}_{t_{j+1}}(b) \Delta_j +  \sum_{j=0}^{\infty} \Delta_j
\widehat{\varphi}_{t_{j+1}}(a) \Delta_j\Delta_{j+1} \widehat{\varphi}_{t_{j+2}}(b)
\Delta_{j+1} \\
&& \ \ \ \ \ \ \ \ \ \ \ \ \ \ \ \ \ \ \ \ \ \ \ \ \ +
\sum_{j=0}^{\infty} \Delta_{j+1} \widehat{\varphi}_{t_{j+2}}(a) \Delta_{j+1}\Delta_j
\widehat{\varphi}_{t_{j+1}}(b) \Delta_j
\ \ \ \ \ \ \ \ \ \ \ \ \ \ \text{(by (\ref{BB}))} \\
&& = \sum_{j=1}^{\infty} \Delta_j \widehat{\varphi}_{t_{j+1}}(a)
\widehat{\varphi}_{t_{j+1}}(b) \Delta_j^2  \Delta_j  +  \sum_{j=1}^{\infty}
\Delta_j \widehat{\varphi}_{t_{j+1}}(a) \widehat{\varphi}_{t_{j+2}}(b)
\Delta_{j+1}\Delta_{j}\Delta_{j+1} \\
&& \ \ \ \ \ \ \ \ \ \ \ \ \ \ \ \ \ \ \ \ \ \ \ \ \
+\sum_{j=1}^{\infty} \Delta_{j} \widehat{\varphi}_{t_{j+1}}(a)
\widehat{\varphi}_{t_{j}}(b) \Delta_{j-1} \Delta_j \Delta_{j-1}
\ \ \ \ \ \ \ \ \ \text{(using (t2) and (p1))} \\
&& = \sum_{j=1}^{\infty} \Delta_j \widehat{\varphi}_{t_{j+1}}(a)
\widehat{\varphi}_{t_{j+1}}(b) \Delta_j^2  \Delta_j  +  \sum_{j=1}^{\infty}
\Delta_j \widehat{\varphi}_{t_{j+1}}(a) \widehat{\varphi}_{t_{j+1}}(b)
\Delta_{j+1}\Delta_{j} \Delta_{j+1}\\
&& \ \ \ \ \ \ \ \ \ \ \ \ \ \ \ \ \ \ \ \ \ \ \ \ \
+\sum_{j=1}^{\infty} \Delta_{j}
\widehat{\varphi}_{t_{j+1}}(a) \widehat{\varphi}_{t_{j+1}}(b)
\Delta_{j-1} \Delta_j  \Delta_{j-1}
\ \ \ \ \ \ \ \ \ \ \ \text{(using (t1))}
\\
&& = \sum_{j=1}^{\infty} \Delta_j \widehat{\varphi}_{t_{j+1}}(ab) \Delta_j^2
\Delta_j  +  \sum_{j=1}^{\infty} \Delta_j \widehat{\varphi}_{t_{j+1}}(ab)
\Delta_{j+1}\Delta_{j} \Delta_{j+1} \\
&& \ \ \ \ \ \ \ \ \ \ \ \ \ \ \ \ \ \ \ \ \ \ \ \ \
+\sum_{j=1}^{\infty} \Delta_{j}
\widehat{\varphi}_{t_{j+1}}(ab) \Delta_{j-1} \Delta_j
\Delta_{j-1}
\ \ \ \ \ \ \ \ \ \ \ \ \ \ \ \ \ \ \
\text{(using (t4))} \\
&& = \sum_{j=0}^{\infty} \Delta_j \widehat{\varphi}_{t_{j+1}}(ab)  \Delta_j
 \ \ \  \ \ \ \ \ \ \ \ \  \ \ \ \ \ \ \text{(since $(\Delta_j^2 +
\Delta_{j+1}^2 + \Delta_{j-1}^2)\Delta_j = \Delta_j$,  $j \geq 2$)}.
\end{eqnarray*}
Hence $\varphi^f(ab) = \varphi^f(a)\varphi^f(b)$ for all $a,b \in D$. The
other conditions for $\varphi^f$ to be a $*$-homomorphism (i.e. linearity
and self-adjointness) are established in the same way.
\end{proof}

\begin{rmk}
{\rm
If $\varphi:D\to Q(E)$ is a genuine $*$-homomorphism then it is natural to
use a $t$-independent lift $\widehat{\varphi}$; then it is obvious that
$\varphi^f=\varphi$, i.e. our basic construction does not change
$*$-homomorphisms.

}
\end{rmk}

We will refer to $\varphi^f$ as \emph{a folding} of $\varphi$. The quadruple
$$
\left( \left(\widehat{\varphi}_t\right)_{t \in [1,\infty)}, \{u_n\}_{n=0}
^{\infty}, \{n_i\}_{i=0}^{\infty}, \{t_i\}_{i=1}^{\infty}\right)
$$
which goes into the construction of $\varphi^f$ will be called \emph{the
folding data}.

At first sight it is not clear how much a folding depends on the folding
data chosen for its construction. Furthermore, it is not difficult to vary
the construction in different ways, for example by omitting condition (t3).
A major part of the proof consists of showing that in the appropriate
setting and modulo the appropriate equivalence relations the construction is
in fact independent of all choices made. Specifically, we will show that it
gives rise to a group homomorphism
$$
F : [[SD,Q(E) \otimes \mathbb K]]  \to [SD, Q(E \otimes \mathbb K)\otimes
\mathbb K].
$$
For this purpose, condition (t3) will come in handy.

\section{Preparing the ground}\label{s4}

Let $A$ and $B$ be $C^*$-algebras, $A$ separable, $B$ $\sigma$-unital. We
say that an asymptotic homomorphism $\varphi = (\varphi_t)_{t \in
[1, \infty)} : A \to B$ is \emph{uniformly continuous} when the function
$[1, \infty) \ni t \mapsto \varphi_t(a)$ is uniformly continuous for all
$a \in A$, i.e. when the following holds:
$$
\forall a \in A \ \forall \epsilon > 0 \ \exists \delta > 0 \ : \ |t-s| \leq
\delta \ \Rightarrow \ \|\varphi_t(a) - \varphi_s(a)\| \leq \epsilon .
$$

\begin{lem}\label{000} Let $\varphi, \psi : A \to B$ be asymptotic
homomorphisms such that $\lim_{t \to \infty} \varphi_t(a) - \psi_t(a) = 0$
for all $a \in A$. If $\psi$ is uniformly continuous then so is $\varphi$.
\end{lem}
\begin{proof} Left to the reader.
\end{proof}

\begin{lem}\label{repar1} Let $\varphi : A \to B$ be an asymptotic
homomorphism. There exists a continuous increasing function $r : [1, \infty)
\to [1, \infty)$ such that $\lim_{t \to \infty}r(t) = \infty$ and
$(\varphi_{r(t)})_{t \in [1, \infty)} : A \to B$ is uniformly continuous.
\end{lem}

\begin{proof} Thanks to Lemma \ref{000} it suffices to prove the statement
when $\varphi$ is equi-continuous.
Let $F_1 \subseteq F_2 \subseteq F_3 \subseteq
\dots $ be a sequence of finite sets with dense union in $A$. There
is a $\delta_n > 0$ so small that $\|\varphi_t(a) - \varphi_s(a)\|
\leq \frac{1}{n}$ for all $a \in F_n$ and all $t,s \in [1, n]$ with
$|t-s| \leq \delta_n$. For each $n$ we choose $k_n \in \mathbb N$
so large that $\delta_{n+4} \geq \frac{1}{k_n}$. Set $N_j = 1 +
\sum_{n=1}^j k_n$, $j = 1,2, \dots $, and $N_0 = 1$. Define
$r : [1, \infty) \to [1, \infty)$ by
$$
r(t) = j + \frac{t- N_{j-1}}{k_j} \ , \ t \in [N_{j-1}, N_j] .
$$
To check that $(\varphi_{r(t)})_{t \in [1, \infty)}$ is uniformly
continuous, let $a \in A$ and $\epsilon >0$. Since $\varphi$ is
equi-continuous there is an $n \in \mathbb N$ and an element $b \in F_n$
such that $\frac{1}{n} \leq \frac{\epsilon}{3}$ and $\sup_{t \in [1, \infty)}
\|\varphi_t(a) -
\varphi_t(b)\| \leq \frac{\epsilon}{3}$. For $x \geq N_n$ and $|x-y| \leq 1$,
we have that $|r(x)-r(y)| \leq \frac{1}{k_j} \leq \delta_{j+4}$ when $x
\in [N_j,N_{j+1}]$. Since $r(x),r(y) \in [1,j+4]$ in this case, we see that
$\|\varphi_{r(x)}(b) - \varphi_{r(y)}(b)\| \leq \frac{1}{n+4} \leq
\frac{\epsilon}{3}$ for all $x,y$ with $x \geq N_n$ and $|x-y| \leq 1$.
Choose $\delta \in ]0,1]$ so small that $\|\varphi_{r(x)}(b) -
\varphi_{r(y)}(b)\| \leq \frac{\epsilon}{3}$ for all $x,y \in [1,N_n + 1]$
with $|x-y| \leq \delta$. Then $\|\varphi_{r(x)}(b) - \varphi_{r(y)}(b)\|
\leq \frac{\epsilon}{3}$ for all $x,y \in [1,\infty)$ with $|x-y| \leq
\delta$. It follows that $\|\varphi_{r(x)}(a) - \varphi_{r(y)}(a)\| \leq
\epsilon$ for all $x,y \in [1,\infty)$ with $|x-y| \leq \delta$.
\end{proof}

\begin{lem}\label{crux2} Let $\lambda : A \to Q(B)$ be an equi-continuous
asymptotic extension. Assume that $\lambda$ is homotopic to a split
asymptotic extension, i.e. that there is a homotopy $\Phi:A\to IQ(B)$
connecting $\lambda=\ev_0\circ\Phi$ with an asymptotic extension
$\ev_1\circ\Phi$ which is split (as an asymptotic extension). It follows
that there is an equi-continuous asymptotic homomorphism $\delta :
A \to M(B)$ such that $\lambda_t = q_B \circ \delta_t$ for all $t$.
\end{lem}

\begin{proof} The proof is based on an idea of Voiculescu, cf. \cite{V}, and
is a refinement of the proof of Lemma 2.1 in \cite{MT3}. By assumption there
is an asymptotic homomorphism $\psi' : A \to M(B)$ such that $\lambda$
is homotopic to the asymptotic extension $q_B \circ \psi'_t,t \in [1,\infty)$.
Let $\psi : A \to M(B)$ be an equi-continuous asymptotic homomorphism such
that $\lim_{t \to \infty} \psi'_t(a) - \psi_t(a) = 0$ for all $a \in A$.
There is then an equi-continuous asymptotic homomorphism $\Phi : A \to
IQ(B)$ which is a homotopy connecting $\lambda$ and $q_B \circ \psi$, i.e.
$\ev_0 \circ \Phi_t = \lambda_t$ and $\ev_1 \circ \Phi_t = q_B \circ
\psi_t$ for all $t$. Define $\Lambda : A \to C_b([1,\infty), IQ(B))$ by
$\Lambda(a)(t) = \Phi_t(a)$. $\Lambda$ is continuous since $\Phi$ is
equi-continuous. Let $Q : C_b([1,\infty), IM(B)) \to C_b([1,\infty),IQ(B))$
be the surjective $*$-homomorphism induced by $\id_{C[0,1]} \otimes q_B :
IM(B) \to IQ(B)$, i.e.
$Q(f)(t) = \left(\id_{C[0,1]} \otimes q_B\right)(f(t))$,
$t\in [1,\infty)$. By the Bartle-Graves selection theorem there is a
continuous section $S$ for $Q$. We set
$$
\mu^s_t(a) = \ev_s\left(S \circ \Lambda(a)(t)\right) + s\left(\psi_t(a) - \ev_1\left(S \circ \Lambda(a)(t)\right)\right)
$$
for all $a \in A$, all $t \in [1,\infty)$ and all $s \in [0,1]$. Then
$$
q_B\left( \mu^0_t(a)\right) = \ev_0\left(Q \circ S \circ \Lambda(a)(t)
\right) = \lambda_t(a)
$$
for all $a,t$. Hence $\mu^s_t, s \in [0,1], t \in [1,\infty)$, is an
equi-continuous family of maps such that
\begin{enumerate}
\item[(m1)] $(\mu_t^1)_{t \in [1,\infty)}$ is an asymptotic homomorphism,
and
\item[(m2)] $q_B \circ \mu_t^0(a) = \lambda_t(a)$ for all $a \in A$ and all
$t$.
\end{enumerate}
Furthermore, since $\psi_t(a) - \ev_1\left(S \circ \Lambda(a)(t)\right) \in B$,
\begin{enumerate}
\item[(m3)] $\lim_{t \to \infty} \sup_{s \in [0,1]} \|q_B \circ \mu^s_t(a^*)
- q_B \circ \mu^s_t(a)^*\| = 0$,
\item[(m4)] $\lim_{t \to \infty} \sup_{s \in [0,1]} \| q_B \circ \mu^s_t(a)
q_B \circ \mu^s_t(b) -  q_B \circ \mu^s_t(ab)\| = 0$, and
\item[(m5)] $\lim_{t \to \infty} \sup_{s \in [0,1]} \| q_B \circ \mu^s_t(a)
+ \lambda q_B \circ \mu^s_t(b) -  q_B \circ \mu^s_t(a +\lambda b)\| = 0$,
\end{enumerate}
for all $a,b \in A$ and all $\lambda \in \mathbb C$.
Choose continuous functions $f_i : [1, \infty) \to [0,1]$,
$i = 0,1,2,\dots$,
such that
\begin{enumerate}
\item[(m6)] $f_0(t) =0$ for all $t \in [1,\infty)$,
\item[(m7)] $f_n \leq f_{n+1}$ for all $n$,
\item[(m8)] for each $ n \in \mathbb N$, there is an $m_n \in
\mathbb N$ such that
$f_i(t) = 1$ for all $i \geq m_n$, and all $t \in [1,n+1]$,
\item[(m9)] $\lim_{t \to \infty} \max_i |\mu_t^{f_i(t)}(a) -
\mu_t^{f_{i+1}(t)}(a)| = 0$ for all $a \in A$.
\end{enumerate}
The last requirement, (m9), can be fulfilled thanks to the separability
of $A$
and the equi-continuity of the family $\mu^s_t$, $s \in [0,1]$, $t \in
[1, \infty)$.
Let $F_1 \subseteq F_2 \subseteq F_3 \subseteq \dots $ be an increasing
sequence of finite subsets with dense union in $A$. For each $n$, choose
$m_n \in \mathbb N$ as in
(m8). We may assume that $m_{n+1} > m_n$. By general facts
on quasi-central approximate units (\cite{A}, \cite{P}) we can choose
elements
$$
X^n_{0}   \geq X^n_{1} \geq X^n_{2} \geq  \dots
$$
in $B$ such that $0 \leq X^n_i \leq 1$ for all $i$ and $X^n_i = 0$ for $i
\geq  m_n$, and
\begin{enumerate}
\item[(i)] $X^n_{i}X^n_{i+1} = X^n_{i+1}$ for all $i$,
\item[(ii)] $\|X^n_i x - x\| \leq \frac{1}{n} + \|q_B(x)\|$ for all $i =
0,1,2, \dots , m_n - 1$, and all $x \in S_n$,
\item[(iii)] $\|X^n_i y - y X^n_i\| \leq \frac{1}{n}$ for all $ i$ and all
$y \in L_n$,
\end{enumerate}
where $L_n = \{ \mu_t^s(a) : \ s \in
[0,1], \ a \in F_n, \ t \in [1,n +1] \}$ and
\begin{eqnarray*}
S_n &=& \{\mu_t^s(a) + \lambda \mu_t^s(b) - \mu_t^s(a + \lambda b) : \ a,b
\in F_n, \ s \in [0,1], \ t \in [1,n+1], \ \lambda \in \mathbb C, |\lambda|
\leq n \} \\
& \cup& \{\mu_t^s(ab) - \mu_t^s(a)\mu_t^s(b) : \ a,b \in F_n, \ s \in [0,1],
 \ t \in [1,n+1] \}  \\
& \cup& \{\mu_t^s(a^*) - {\mu_t^s(a)}^* : \ a \in F_n, \ s \in [0,1], \ t \in
[1, n+1] \}
\end{eqnarray*}
are the compact subsets of $M(B)$.

By choosing the $X^n_i$'s recursively, we can arrange that $X^{n+1}_i X^
n_k = X^n_k$ for all $k$ and all $i$. By connecting first $X^n_0$ to $X^{
n+1}_0$ via the straight line between them, then $X^n_{1}$ to $X^{n+1}_{1}$
via a straight line, then $X^n_{2}$ to $X^{n+1}_{2}$ etc., we obtain
norm-continuous paths,
$X(t,i)$, $t \in [n,n+1]$, $i = 0,1,2,3, \dots$, in $B$
such that $X(n,i) = X^n_i, \ X(n+1,i) = X^{n+1}_i$ for all $i$ and
\begin{enumerate}
\item[(m10)] $X(t,{i})X(t,i+1) = X(t,i+1), \  t \in [n,n+1]$, for all $i$,
\item[(m11)] $\|X(t,i) x - x\| \leq \frac{1}{n} + \|q_B(x)\|$ for all $i =
0,1,2, \dots , m_n - 1, \ t \in [n,n+1]$ and all $x \in S_n$,
\item[(m12)] $\|X(t,i) y - y X(t,i)\| \leq \frac{1}{n}$ for all $i$, all $t
\in [n,n+1]$ and all $y \in L_n$.
\end{enumerate}
In addition, $X(t,i) = 0, i \geq m_{n+1}, t \in [n,n+1]$. Set $\Lambda^t_
0 = \sqrt{1 - X(t,0)}$, and $\Lambda^t_j = \sqrt{X(t,j-1)-X(t,j)}, j \geq 1$,
and define $\delta_t : A \to M(B)$ by
$$
\delta_t(a) = \sum_{j = 0}^{\infty} \Lambda^t_j \mu_t^{f_j(t)}(a)
\Lambda^t_j .
$$
Note that the sum is finite for $t$ in a compact set, and that $t \mapsto
\delta_t(a)$ is norm-continuous. Observe also that
\begin{enumerate}
\item[(m13)] $\Lambda^t_j\Lambda^t_i = 0, |i-j| \geq 2$, for all $t$, and
\item[(m14)] $\sum_{j=0}^{\infty} \left(\Lambda^t_j\right)^2 = 1$ for all
$t$.
\end{enumerate}
It follows from (m6) that $q_B \circ \delta_t(a) = q_B \circ \mu^0_t(a)$
for all $a,t$. Thanks to (m2) it now only remains to show that $\delta =
\left( \delta_t\right)_{t \in [1, \infty)}$ is an asymptotic homomorphism.
We check that it is multiplicative. Because
$(\delta_t)_{t\in[1,\infty)}$ is an
equi-continuous family, it suffices to consider $a,b \in F_n$, and show that
$\lim_{t \to \infty} \delta_t(a)\delta_t(b) - \delta_t(ab) = 0$.
For this purpose observe that for any sequence of functions $g_i :
[1,\infty) \to M(B), i = 0,1,2, \dots$, we have the estimates
\begin{equation}\label{ebbe}
\| \sum_{j=0}^{\infty} \Lambda^t_j g_j(t)\Lambda^t_j\| \leq \sup_i \|g_i(t)\|,
\ \ \| \sum_{j=0}^{\infty} \Lambda^t_j g_j(t)\Lambda^t_{j+1}\| \leq
\sup_i \|g_i(t)\|
\end{equation}
for all $t$. This follows as in the proof of Lemma \ref{E1}.

In the following we will write $\sim$ between two expressions that depend
on $t$ when their difference tends to $0$ in norm as $t$ tends to infinity.
We will use repeatedly the estimate (\ref{ebbe})
in the calculation below.
\begin{eqnarray*}
&& \Bigl( \sum_{j=0}^{\infty} \Lambda^t_j \mu_t^{f_j(t)}(a)\Lambda^t_j
\Bigr) \Bigl( \sum_{j=0}^{\infty} \Lambda^t_j \mu_t^{f_j(t)}(b)\Lambda^t_j
\Bigr) \\
&&= \sum_{j=0}^{\infty} \Lambda^t_j \mu_t^{f_j(t)}(a)\left(\Lambda^t_j
\right)^2 \mu_t^{f_j(t)}(b)\Lambda^t_j +  \sum_{j=0}^{\infty} \Lambda^t_j
\mu_t^{f_j(t)}(a){\Lambda^t_j}{\Lambda^t_{j+1}} \mu_t^{f_{j+1}(t)}(b)
\Lambda^t_{j+1} \\
&& \ \ \ \  \ \ \ \ \ \ \ \  \ \ \  \ \  \
+\sum_{j=0}^{\infty} \Lambda^
t_{j+1} \mu_t^{f_{j+1}(t)}(a){\Lambda^t_{j+1}}{\Lambda^t_{j}}
\mu_t^{f_{j}(t)}(b)\Lambda^t_{j}
\ \ \ \ \ \ \  \ \ \ \ \ \text{(using (m13))} \\
&& \sim  \sum_{j=0}^{\infty} \Lambda^t_j \mu_t^{f_j(t)}(a)
\mu_t^{f_j(t)}(b)\left(\Lambda^t_j\right)^2 \Lambda^t_j + \sum_{j=0}^{\infty}
\Lambda^t_j \mu_t^{f_j(t)}(a)\mu_t^{f_{j+1}(t)}(b){\Lambda^t_j}
{\Lambda^t_{j+1}}\Lambda^t_{j+1} \\
&& \ \ \ \  \ \ \ \ \ \ \ \  \ \ \  \ \  \
+\sum_{j=1}^{\infty} \Lambda^t_{j} \mu_t^{f_{j}(t)}(a)\mu_t^{f_{j-1}(t)}(b)
{\Lambda^t_{j}}{\Lambda^t_{j-1}}\Lambda^t_{j-1}
\ \ \  \ \ \ \ \ \ \ \
\text{(using (m12))} \\
&& \sim  \sum_{j=0}^{\infty} \Lambda^t_j \mu_t^{f_j(t)}(a)
\mu_t^{f_j(t)}(b) \left(\Lambda^t_j\right)^2 \Lambda^t_j +
\sum_{j=0}^{\infty} \Lambda^t_j
\mu_t^{f_j(t)}(a)\mu_t^{f_{j}(t)}(b){\Lambda^t_j}{\Lambda^t_{j+1}}
\Lambda^t_{j+1} \\
&& \ \ \ \  \ \ \ \ \ \ \ \  \ \ \  \ \  \
+\sum_{j=
1}^{\infty} \Lambda^t_{j} \mu_t^{f_{j}(t)}(a)\mu_t^{f_{j}(t)}(b){\Lambda
^t_{j}}{\Lambda^t_{j-1}}\Lambda^t_{j-1}
\ \ \ \ \ \ \ \ \  \ \ \ \ \ \text{(using (m9))}
\\
&& \sim \sum_{j=0}^{\infty} \Lambda^t_j \mu_t^{f_j(t)}(ab)\left(\Lambda^t_j
\right)^2 \Lambda^t_j + \sum_{j=0}^{\infty} \Lambda^t_j \mu_t^{f_j(t)}(ab)
{\Lambda^t_j}{\Lambda^t_{j+1}}\Lambda^t_{j+1} \\
&& \ \ \ \  \ \ \ \ \ \ \ \  \ \ \  \ \  \
+  \sum_{j=
1}^{\infty} \Lambda^t_{j} \mu_t^{f_{j}(t)}(ab){\Lambda^t_{j}}
{\Lambda^t_{j-1}}\Lambda^t_{j-1} \ \ \  \ \ \ \ \ \ \ \ \ \
\text{(using (m11),(m4) and (m1))} \\
&& = \sum_{j=0}^{\infty} \Lambda^t_j \mu_t^{f_j(t)}(ab)\Lambda^t_j \ \ \
\ \ \ \ \ \ \ \ \ \ \ \ \ \ \ \ \ \ \ \ \ \ \ \ \ \ \ \ \ \ \
\ \ \ \ \ \ \ \ \ \ \
\text{(using (m13) and (m14)).}
\end{eqnarray*}
Similar considerations show that $\delta$ is also asymptotically linear and
asymptotically self-adjoint.
\end{proof}

The following lemma is crucial in understanding to what extend the basic
construction depends on a choice of folding data.

\begin{lem}\label{XX} Let $\varphi, \psi : A \to Q(B)$ be equi-continuous
asymptotic extensions. Assume that the asymptotic extension
$$
\left( \begin{smallmatrix} \varphi & 0 \\ 0 & \psi \end{smallmatrix} \right
) : A \to M_2(Q(B))
$$
is split (as an asymptotic extension).
It follows that for any foldings, $\varphi^f$ and $\psi^f$, of
$\varphi$ and $\psi$, respectively, the extension
$$
\left( \begin{smallmatrix} \varphi^f & 0 \\ 0 & \psi^f \end{smallmatrix}
\right) : A \to M_2(Q(B))
$$
is asymptotically split.
\end{lem}
\begin{proof} Let $\left( \left(\widehat{\varphi}_t\right)_{t \in [1,\infty)},
\{u_n\}_{n=0}^{\infty}, \{n_i\}_{i=0}^{\infty}, \{t_i\}_{i=1}^{\infty}
\right)$ be the folding data used to define $\varphi^f$ and $\left(
(\widehat{\psi}_t)_{t \in [1,\infty)}, \{u'_n\}_{n = 0}^{\infty},
\{n'_i\}_{i=0}^{\infty}, \{t'_i\}_{i=1}^{\infty}\right)$ the folding data
used to define $\psi^f$.  To simplify notation, set $w_i = u_{n_i}$ and
$w'_i = u'_{n'_i}$ for all $i$. By assumption there is a an asymptotic
homomorphism $\pi=(\pi_t)_{t\in[1,\infty}:A\to M_2(M(B))$ such that
\begin{equation}\label{EF4}
\lim_{t \to \infty} q_{M_2(B)} \circ \pi_t(a) - \left( \begin{smallmatrix}
\varphi_t(a) & \\ & \psi_t(a) \end{smallmatrix} \right) = 0
\end{equation}
for all $a \in A$. We may assume that $\pi$ is an equi-continuous asymptotic
homomorphism. Let $F_1 \subseteq F_2 \subseteq F_3 \subseteq \dots$ be a
sequence of finite sets with dense union in $A$. Since
$\left( \left(\begin{smallmatrix} u_n & \\  & u'_n \end{smallmatrix}
\right) \right)_{n \in \mathbb N}$ is an approximate unit in $M_2(B)$ we can
choose a sequence
$0=r_0 < r_1 < r_2 < r_3 < \dots $ in  $\mathbb N$ such that
\begin{equation}\label{A}
\|\left ( \begin{smallmatrix} 1- u_j & {} \\ {} & 1- u'_j \end{smallmatrix}
\right )\left[ \left ( \begin{smallmatrix} \widehat{\varphi}_{t}(d) & {} \\ {}
& \widehat{\psi}_{t}(d) \end{smallmatrix} \right ) - \pi_{t}(d)\right]\| \leq
\frac{1}{n} + \| \left ( \begin{smallmatrix} \varphi_{t}(d) & {} \\ {} &
\psi_{t}(d) \end{smallmatrix} \right ) - q_{M_2(B)} \circ \pi_{t}(d)\|
\end{equation}
for all $t \in [1,n+2]$, $j \geq r_{n-2}$, $d \in F_n$. It follows from
conditions (p1) and (t4)-(t6) on the folding data that we can arrange that
\begin{enumerate}
\item[(a)]
$\|u_j \widehat{\varphi}_t(a) - \widehat{\varphi}_t(a)u_j\| \leq \frac{1}
{i+1}, \ t \in [1,i+2], \ j \geq r_i$,
\item[(b)]
$\|\left(1-u_j\right) \left( \widehat{\varphi}_t(a) \widehat{\varphi}_t(b)
- \widehat{\varphi}_t(ab)\right)\| \leq \|\varphi_t(a)\varphi_t(b) - \varphi
_t(ab)\| + \frac{1}{i+1}, \ t \in [1,i+2], \ j \geq r_i$,
\item[(c)]
$\|\left(1 - u_j\right) \left( \widehat{\varphi}_t(a) + \lambda \widehat{
\varphi}_t(b) -  \widehat{\varphi}_t(a + \lambda b)\right)\| \leq \|\varphi_t(a
) + \lambda \varphi_t(b) - \varphi_t(a + \lambda b)\| + \frac{1}{i+1}, \ t
\in [1,i+2], \ j \geq r_i$,
\item[(d)]
$\|\left( 1-u_j\right) \left( \widehat{\varphi}_t(a^*) - \widehat{\varphi
}_t(a)^*\right)\| \leq \|\varphi_t(a^*) - \varphi_t(a)^*\| + \frac{1}{i+1},
 \ t \in [1, i+2], \ j \geq r_i$
\end{enumerate}
and
\begin{enumerate}
\item[(a')]
$\|u'_j \widehat{\psi}_t(a) - \widehat{\psi}_t(a)u'_j\| \leq \frac{1}{i+
1}, \ t \in [1,i+2], \ j \geq r_i$,
\item[(b')]
$\|\left(1-u'_j\right) \left( \widehat{\psi}_t(a) \widehat{\psi}_t(b) -
 \widehat{\psi}_t(ab)\right)\| \leq \|\psi_t(a)\psi_t(b) - \psi_t(ab)\| + \frac
{1}{i+1}, \ t \in [1,i+2],\ j \geq r_i$,
\item[(c')]
$\|\left(1 - u'_j\right) \left( \widehat{\psi}_t(a) + \lambda
\widehat{\psi}_t(b) -  \widehat{\psi}_t(a + \lambda b)\right)\| \leq \|\psi_t(a) +
\lambda \psi_t(b) - \psi_t(a + \lambda b)\| + \frac{1}{i+1}, \ t \in [1,i+2],
\ j \geq r_i$,
\item[(d')]
$\|\left( 1-u'_j\right) \left( \widehat{\psi}_t(a^*) - \widehat{\psi}_t(
a)^*\right)\| \leq \|\psi_t(a^*) - \psi_t(a)^*\| + \frac{1}{i+1}, \ t \in [
1, i+2], \ j \geq r_i$,
\end{enumerate}
for all $a,b \in F_i$ and all $\lambda \in \mathbb C, |\lambda| \leq i$.

We will construct continuous paths $w_i(t)$, $w'_i(t)$,
$t \in[0,\infty)$, $i = 0,1,2, \dots $, in $B$ such that
\begin{enumerate}
\item[(w1)]
$\{w_i(t)\}_{i =0}^
{\infty}$ and $\{w'_i(t)\}_{i =0}^{\infty}$ are unit sequences for all $t$,
\item[(w2)]
$w_i(0) = w_i$, $w'_i(0) = w'_i$ for all $i$,
\item[(w3)]
$w_i(t) = w_i$, $w'_i(t) =
w'_i$ for $i \geq r_{k+1}$ when $t \in [k,k+1]$,
\item[(w4)]
$w_0(k) \geq
w_{r_{k-1}}$, $w'_0(k) \geq w'_{r_{k-1}}$,
\item[(w5)]
$w_i(t) \in \co \{ u_j : j \geq n_i\}, \ w'_i(t) \in \co \{ u'_j : j \geq
n' _i\}$ for all $i,t$,
\item[(w6)]
$w_i(t) \in \co \{ u_j : j \geq r_{k-1}\}, \ w'(t) \in \co \{u'_j : j \geq
r_{k-1}\}$ for all $i$ when $t \in [k,k+1]$.
\end{enumerate}
To see how to do that, assume that we
have constructed $\{w_i(t\}_{i=0}^{\infty}$ and $\{w'_i(t)\}_{i=0}^{\infty
}$ for all $t \in [0,k]$. Then $w_{r_{k-1}} \leq w_0(k) \leq w_1(k) \leq
\dots \leq w_{r_{k} -1}(k) \leq w_{r_{k}}(k) = w_{r_{k}}$, and $w_i(k) =
w_i$, $i \geq r_{k}$. Thanks to condition (t3) on the folding data there are
elements $v_0 < v_1 < \dots < v_{r_{k}}$ from $\{u_j :
n_{r_{k}}<j<n_{r_{k}+1} \}$ such that
$$
w_0(k) \leq w_1(k) \leq \dots \leq w_{r_{k}-1}(k) \leq w_{r_k} \leq v_0
\leq v_1 \leq \dots \leq v_{r_{k}} \leq w_{r_{k} +1} \leq w_{r_{k}+2}  \leq
\dots
$$
is a unit sequence. We define $w_i(t)$, $t \in [k,k+1]$, such that $w_i(t)
= w_i$, $i \geq r_{k}+1$. To define $w_i(t)$ for $i \leq r_k$ and $t \in [k
,k+1]$, set $I_j = [k + \frac{j}{r_k+1}, k + \frac{j+1}{r_k+1}]$, $j = 0,
1, \dots, r_k$. On the interval $I_j$, $w_{r_k -j}$ is the line from $w_{r_
k-j}(k)$ to $v_{r_k-j}$, while $w_m$ is constant on $I_j$ for $m \neq r_k-j
$. This completes the construction of $\{w_i(t)\}_{i =0}^{\infty}$
for $t \in [k,k+1]$.
Proceed inductively to construct $\{w_i(t)\}_{i =0}^{\infty}$, $t
\in [1,\infty)$.  $\{w'_i(t)\}_{i =0}^{\infty}$, $t \in [1,\infty)$, is
constructed in the same way.

Note that (a)-(d) and (a')-(d') in combination with (w6) imply that
\begin{enumerate}
\item[(A)]
$ \sup_j \|w_j(s) \widehat{\varphi}_t(a) - \widehat{\varphi}_t(a)w_j(s)\|
 \leq \frac{1}{i}, \ s,t \in [i,i+1]$,
\item[(B)]
$\sup_j \|\left(1-w_j(s)\right) \left( \widehat{\varphi}_t(a)
\widehat{\varphi}_t(b) -  \widehat{\varphi}_t(ab)\right)\| \leq \|\varphi_t(a)\varphi_t(b)
- \varphi_t(ab)\| + \frac{1}{i}, \ s,t \in [i,i+1]$,
\item[(C)]
$\sup_j \|\left(1 - w_j(s)\right) \left( \widehat{\varphi}_t(a) +
\lambda \widehat{\varphi}_t(b) -  \widehat{\varphi}_t(a + \lambda b)\right)\| \leq
\|\varphi_t(a) + \lambda \varphi_t(b) - \varphi_t(a + \lambda b)\| + \frac{1}{
i}, \ s,t \in [i,i+1], \ j \geq r_i$,
\item[(D)]
$\sup_j\|\left( 1-w_j(s)\right) \left( \widehat{\varphi}_t(a^*) -
\widehat{\varphi}_t(a)^*\right)\| \leq \|\varphi_t(a^*) - \varphi_t(a)^*\| + \frac
{1}{i}, \ s,t \in [i, i+1], \ j \geq r_i$
\end{enumerate}
and
\begin{enumerate}
\item[(A')]
$\sup_j\|w'_j(s) \widehat{\psi}_t(a) - \widehat{\psi}_t(a)w'_j(s)\| \leq
 \frac{1}{i}, \ s,t \in [i,i+1]$,
\item[(B')]
$\sup_j\|\left(1-w'_j(s)\right) \left( \widehat{\psi}_t(a) \widehat{\psi
}_t(b) -  \widehat{\psi}_t(ab)\right)\| \leq \|\psi_t(a)\psi_t(b) - \psi_t(ab)
\| + \frac{1}{i}, \ s,t \in [i,i+1]$,
\item[(C')]
$\sup_j\|\left(1 - w'_j(s)\right) \left( \widehat{\psi}_t(a) +
\lambda \widehat{\psi}_t(b) -  \widehat{\psi}_t(a + \lambda b)\right)\| \leq \|\psi_t(a)
+ \lambda \psi_t(b) - \psi_t(a + \lambda b)\| + \frac{1}{i}, \ s,t \in [i,
i+1]$,
\item[(D')]
$\sup_j\|\left( 1-w'_j(s)\right) \left( \widehat{\psi}_t(a^*) - \widehat
{\psi}_t(a)^*\right)\| \leq \|\psi_t(a^*) - \psi_t(a)^*\| + \frac{1}{i}, \
s,t \in [i, i+1]$,
\end{enumerate}
for all $a,b \in F_i$ and all $\lambda \in \mathbb C, |\lambda| \leq i$.

Set $\Delta_0(t) = \sqrt{w_0(t)}$, $\Delta'_0(t) = \sqrt{w'_0(t)}$,
$\Delta_i(t) = \sqrt{w_i(t) - w_{i-1}(t)}$,
$\Delta'_i(t) = \sqrt{w'_i(t) - w'_
{i-1}(t)}$, $i \geq 1$, and define $f^i,g^i : [1,\infty) \to [1,\infty)$
by $f^i(t) = \max\{t, t_{i+1}\}$ and $g^i(t) = \max \{t, t'_{i+1}\}$, $i =
1,2, \dots $. Let
$$
\pi_t = \left ( \begin{matrix} \pi_t^{11} &  \pi_t^{12} \\ \pi_t^{21} &
\pi_t^{22} \end{matrix} \right )
$$
be the matrix decomposition of $\pi_t$. We define
$\rho=(\rho_t)_{t\in[1,\infty)}:A \to M_2(M(B))$ by
\begin{eqnarray*}
\rho_t(a) &=& \left( \begin{matrix} \sum_{i=1}^{\infty} \Delta_i(t) \widehat
{\varphi}_{f^i(t)}(a) \Delta_i(t)  & 0  \\  0  & \sum_{i=1}^{\infty}
\Delta'_i(t) \widehat{\varphi}_{g^i(t)}(a) \Delta'_i(t)  \end{matrix} \right) \\
& 
+& \left( \begin{matrix} \Delta_0(t) \pi_t^{11}(a)\Delta_0(t) &  \Delta_0(
t) \pi_t^{12}(a)\Delta'_0(t) \\   \Delta'_0(t) \pi_t^{21}(a)\Delta_0(t)   &
\Delta'_0(t) \pi_t^{22}(a)\Delta'_0(t)  \end{matrix} \right) .
\end{eqnarray*}
Note first of all that due to (w3) one has
\begin{equation}\label{Uii1}
q_{M_2(B)} \circ \rho_t = \left( \begin{smallmatrix} \varphi^f & 0 \\ 0 &
\psi^f \end{smallmatrix} \right)
\end{equation}
for all $t$. It suffices therefore to show that $\rho$ is an asymptotic
homomorphism.
To simplify the following calculations set $\Phi_t(a) = \sum_{i=1}^{
\infty} \Delta_i(t) \widehat{\varphi}_{f^i(t)}(a) \Delta_i(t)$, $\Psi_t(a) =
\sum_{i=1}^{\infty} \Delta'_i(t) \widehat{\psi}_{g^i(t)}(a) \Delta'_i(t)$ so that
$$
\rho_t(a) = \left( \begin{smallmatrix} \Phi_t(a) & 0 \\ 0 & \Psi_t(a)
\end{smallmatrix} \right) + \left( \begin{smallmatrix} \Delta_0(t) & 0 \\ 0 &
\Delta'_0(t) \end{smallmatrix} \right)\pi_t(a) \left( \begin{smallmatrix}
\Delta_0(t) & 0 \\ 0 & \Delta'_0(t) \end{smallmatrix} \right) .
$$
Note that $(\Phi_t)_{t \in [1,\infty)}$ is an equi-continuous
family by Lemma
\ref{E1} since $(\widehat{\varphi}_t)_{t\in [1,\infty)}$ is.
As $\pi$ is also
equi-continuous, we find that $(\rho_t)_{t \in [1,\infty)}$ is an
equi-continuous family. To show that $\rho_t(a) \rho_t(b) \sim \rho_t(ab)$
for all $a,b$ it suffices therefore to check for $a,b \in F_n$.
We find that
\begin{eqnarray*}
&&\Phi_t(a) \Phi_t(b)  \\
&& = \sum_{i=1}^{\infty} \Delta_i(t)\widehat{\varphi}_{f^i(t)}(a) \Delta_i(t
)^2 \widehat{\varphi}_{f^i(t)}(b)\Delta_i(t) +  \sum_{i=1}^{\infty} \Delta_i(
t)\widehat{\varphi}_{f^i(t)}(a) \Delta_i(t)\Delta_{i+1}(t) \widehat{\varphi}_{f^{i+
1}(t)}(b)\Delta_{i+1}(t)\\
&& \ \ \ \ \ \ \ \ \ \ \ \ \ \ + \sum_{i=2}^{\infty} \Delta_i(t)\widehat{\varphi}
_{f^i(t)}(a) \Delta_i(t)\Delta_{i-1}(t) \widehat{\varphi}_{f^{i-1}(t)}(b)
\Delta_{i-1}(t)\\
&& \sim  \sum_{i=1}^{\infty} \Delta_i(t)\widehat{\varphi}_{f^i(t)}(a)
\widehat{\varphi}_{f^i(t)}(b)\Delta_i(t)^2\Delta_i(t) +  \sum_{i=1}^{\infty}
\Delta_i(
t)\widehat{\varphi}_{f^i(t)}(a) \widehat{\varphi}_{f^{i+1}(t)}(b) \Delta_i(t)\Delta
_{i+1}(t)^2 \\
&& \ \ \ \ \ \ \ \ \ \ \ \ \ \ + \sum_{i=2}^{\infty} \Delta_i(t)
\widehat{\varphi}_{f^i(t)}(a) \widehat{\varphi}_{f^{i-1}(t)}(b)
\Delta_i(t)\Delta_{i-1}(t)^2
\ \ \ \ \ \ \ \ \ \text{(by (A), (w5) and (p1))}\\
&& \sim  \sum_{i=1}^{\infty} \Delta_i(t)\widehat{\varphi}_{f^i(t)}(a)
\widehat{\varphi}_{f^i(t)}(b)\Delta_i(t)^2\Delta_i(t) +  \sum_{i=1}^{\infty}
\Delta_i(
t)\widehat{\varphi}_{f^i(t)}(a) \widehat{\varphi}_{f^{i}(t)}(b) \Delta_i(t)\Delta_{
i+1}(t)^2 \\
&& \ \ \ \ \ \ \ \ \ \ \ \ \ \ + \sum_{i=2}^{\infty} \Delta_i(t)
\widehat{\varphi}_{f^{i}(t)}(a) \widehat{\varphi}_{f^{i}(t)}(b)
\Delta_i(t)\Delta_{i-1}(t)^2
\ \ \ \ \ \ \ \ \ \ \ \ \ \ \ \ \ \ \text{(by (t1))}\\
&& \sim  \sum_{i=1}^{\infty} \Delta_i(t)\widehat{\varphi}_{f^i(t)}(ab)\Delta_i
(t)^2\Delta_i(t) +  \sum_{i=1}^{\infty} \Delta_i(t)\widehat{\varphi}_{f^i(t)}
(ab) \Delta_i(t)\Delta_{i+1}(t)^2 \\
&& \ \ \ \ \ \ \ \ \ \ \ \ \ \ + \sum_{i=2}^{\infty}
\Delta_i(t)\widehat{\varphi}_{f^{i}(t)}(ab)\Delta_i(t)\Delta_{i-1}(t)^2
\ \ \ \ \ \ \ \ \ \ \ \ \ \ \ \ \text{(by (D), (w5) and (t4))}\\
&& = \Phi_t(ab) - \Delta_1(t)^2 \Delta_0(t)^2\widehat{\varphi}_{f^1(t)}(ab) \
\ \ \ \ \ \ \ \ \ \ \ \ \ \ \ \ \ \ \ \ \ \ \ \ \ \ \ \ \
\text{(since $\sum_{j=0}^{\infty}
\Delta_j(t)^2 =1$).}
\end{eqnarray*}
In the same way we see that
$$
\Psi_t(a) \Psi_t(b) \sim \Psi_t(ab) -  \Delta'_1(t)^2 \Delta'_0(t)^2
\widehat{\psi}_{g^1(t)}(ab).
$$
It follows from (\ref{A}) and (w6) that
\begin{equation}\label{poul}
\left( \begin{smallmatrix} \Delta_1(t)^2 & \\  & \Delta'_1(t)^2
\end{smallmatrix} \right)  \left( \begin{smallmatrix}
\widehat{\varphi}_{f^1(t)}(ab) & \\
 & \widehat{\psi}_{g^1(t)}(ab) \end{smallmatrix} \right) \sim \left(
\begin{smallmatrix} \Delta_1(t)^2 & \\  & \Delta'_1(t)^2 \end{smallmatrix}
\right) \pi_t(ab).
\end{equation}
Thus
\begin{equation}\label{Ruth}
\left( \begin{smallmatrix} \Phi_t(a) & \\  & \Psi_t(a) \end{smallmatrix}
\right) \left( \begin{smallmatrix} \Phi_t(b) & \\  & \Psi_t(b)
\end{smallmatrix} \right) \sim \left( \begin{smallmatrix} \Phi_t(ab) & \\
& \Psi_t(ab) \end{smallmatrix} \right) - \left( \begin{smallmatrix}
\Delta_1(t)^2\Delta_0(t)^2 & \\  & \Delta'_1(t)^2\Delta'_0(t)^2
\end{smallmatrix} \right) \pi_t(ab).
\end{equation}
It follows from (A), (A') and (w6) that $ \left( \begin{smallmatrix}
\Delta_0(t) & 0 \\ 0 & \Delta'_0(t) \end{smallmatrix} \right)$ asymptotically
commutes with $\pi_t(a)$ for all $a \in  \bigcup_n F_n$. Hence
\begin{equation}\label{niver}
\left( \begin{smallmatrix} \Delta_0(t) & 0 \\ 0 & \Delta'_0(t)
\end{smallmatrix} \right)\pi_t(a) \left( \begin{smallmatrix} \Delta_0(t) & 0 \\ 0 &
\Delta'_0(t) \end{smallmatrix} \right) \left( \begin{smallmatrix} \Delta_0(t)
& 0 \\ 0 & \Delta'_0(t) \end{smallmatrix} \right)\pi_t(b)
\left( \begin{smallmatrix} \Delta_0(t) & 0 \\ 0 & \Delta'_0(t) \end{smallmatrix}
\right) \sim  \left( \begin{smallmatrix} \Delta_0(t)^2 & 0 \\ 0 & \Delta'_0(t)^2
\end{smallmatrix} \right)\pi_t(ab) .
\end{equation}
Finally, it follows easily from (\ref{poul}) that
\begin{equation}\label{kliver}
\left( \begin{smallmatrix} \Phi_t(a) & \\  & \Psi_t(a) \end{smallmatrix}
\right)  \left( \begin{smallmatrix} \Delta_0(t)^2 & 0 \\ 0 & \Delta'_0(t)^2
\end{smallmatrix} \right) \pi_t(b) \sim  \left( \begin{smallmatrix} \Delta_0
(t)^2 \Delta_1(t)^2 & 0 \\ 0 & \Delta'_0(t)^2\Delta'_1(t)^2 \end{smallmatrix}
\right)\pi_t(ab)
\end{equation}
and
\begin{equation}\label{hans}
  \left( \begin{smallmatrix} \Delta_0(t)^2 & 0 \\ 0 & \Delta'_0(t)^2
\end{smallmatrix} \right) \pi_t(a) \left( \begin{smallmatrix} \Phi_t(b) & \\  &
\Psi_t(b) \end{smallmatrix} \right)  \sim  \left( \begin{smallmatrix} \Delta
_0(t)^2 \Delta_1(t)^2 & 0 \\ 0 & \Delta'_0(t)^2\Delta'_1(t)^2
\end{smallmatrix} \right)\pi_t(ab).
\end{equation}
Now the desired conclusion, that $\rho_t(a)\rho_t(b) \sim \rho_t(ab)$,
follows by combining (\ref{Ruth}), (\ref{niver}), (\ref{kliver}) and
(\ref{hans}).  That $\rho$ is also asymptotically linear and self-adjoint
follows fairly straightforwardly from (C), (C'), (D), (D'), (t5), (t6),
(w6) and the fact that $\pi$ is asymptotically linear and
self-adjoint.
\end{proof}

\begin{lem}\label{unitary} Let $\varphi, \psi : A \to Q(B)$ be asymptotic
extensions that are
unitarily equivalent. It follows that there are foldings,
$\varphi^f$ and $\psi^f$, of $\varphi$ and $\psi$, respectively, such that
$\left( \begin{smallmatrix} \varphi^f & 0 \\ 0 & 0 \end{smallmatrix}
\right)$ and $\left( \begin{smallmatrix} \psi^f & 0 \\ 0 & 0 \end{smallmatrix}
\right)$ are unitarily equivalent.
\end{lem}
\begin{proof} By assumption there is a norm-continuous path
$(v_t)_{t \in [1,\infty)}$,
of unitaries in $M(B)$ such that $\lim_{t \to \infty} \Ad q_B(v_t
) \circ \varphi_t(a) - \psi_t(a) = 0$ for all $a \in A$. It is easy to see
that we can take folding data $\Bigl( (\widehat{\varphi}_t)_{t
\in [1,\infty)}, \{u_n\}_{n=0}^{\infty}, \{n_i\}_{i=0}^{\infty},
\{t_i\}_{i=1}^{\infty} \Bigr)$ for $\varphi$ and
$\Bigl( (\widehat{\psi}_t)
_{t \in [1,\infty)}, \{u'_n\}_{n=0}^{\infty}, \{n'_i\}_{i=0}^{\infty
}, \{t'_i\}_{i=1}^{\infty}\Bigr)$ for $\psi$ such that $u_n = u'_n, t_
n= t'_n$ for all $n$ and $n_i = n'_i$ for all $i$, and such that,
furthermore, $\lim_{i \to \infty} v_{t_{i+1}} - v_{t_i} = 0$ and
$\lim_{i \to \infty} u_{n_i}v_{t_i} - v_{t_i}u_{n_i} = 0$. Let $\varphi^f$
and $\psi^f$ be
the foldings obtained with such choices. Set $\Delta_0 = \sqrt{u_{n_0}}
$ and $\Delta_j = \sqrt{u_{n_j} - u_{n_{j-1}}}, \ j \geq 1$. Then $\sum_{
j=0}^{\infty} \Delta_j v_{t_{j+1}}\Delta_j$ converges in the strict topology
to an element $v \in M(B)$. Using Lemma \ref{E1} we see that $q_B(v)$ is a
unitary in $Q(B)$ such that $\Ad q_B(v) \circ \varphi^f = \psi^f$. Let
$w \in M_2(M(B))$ be a unitary lift of $\left( \begin{smallmatrix} q_B(v)
 & 0 \\ 0 & q_B(v)^* \end{smallmatrix} \right)$. Then
$$
\Ad q_{M_2(B)}(w) \circ \left( \begin{smallmatrix} \varphi^f & 0 \\ 0 & 0
\end{smallmatrix} \right) = \left( \begin{smallmatrix} \psi^f & 0 \\ 0 & 0
\end{smallmatrix} \right).
$$
\end{proof}

We assume now that $B$ is stable. This gives space to add extensions and
asymptotic extensions: Let $W_1,W_2$ be isometries in $M(B)$ such that $W_1W_
1^* + W_2W_2^* = 1$. When $\psi, \varphi : A \to Q(B)$ are asymptotic
extensions we set $(\psi \oplus \varphi)_t(a) = q_B(W_1) \psi_t(a)q_B(W_1^*)
+ q_B(W_2)\varphi_t(a)q_B(W_2^*)$. Up to unitary equivalence this addition
is independent of the choice of isometries, $W_1,W_2$.

\begin{lem}\label{cp} Let $\psi : A \to Q(B)$ be a completely positive
asymptotic extension. Assume that $\psi$ is homotopic to an
asymptotic extension which is split (as an asymptotic extension).
It follows that there are
\begin{enumerate}
\item[1)] a completely positive asymptotic extension $\mu : A \to Q(B)$
which is split (as an asymptotic extension),
\item[2)] a continuous increasing function $r : [1, \infty) \to [1,\infty)$
 with $\lim_{t \to \infty} r(t) = \infty$ , and
\item[3)] a norm-continuous path $(U_t)_{t \in [1, \infty)}$, of unitaries in
$M_2(M(B))$
\end{enumerate}
such that
$$
\lim_{t \to \infty} \Ad q_{M_2(B)}\left(U_t\right) \circ \left(
\begin{smallmatrix} \psi_{r(t)}(a) \\  & \mu_t(a) \end{smallmatrix}
\right) -  \left( \begin{smallmatrix} 0
\\  & \mu_t(a) \end{smallmatrix} \right) = 0
$$
for all $a \in A$.
\end{lem}
\begin{proof} Since $\psi$ is equi-continuous, Lemma \ref{crux2} tells us
that there is an equi-continuous asymptotic homomorphism $\varphi : A \to M(B
)$ such that $q_B \circ \varphi_t = \psi_t$ for all $t \in [1,\infty)$. By
Lemma \ref{repar1} there is a continuous increasing function $r_0 : [1,
\infty) \to [1,\infty)$ such that $\lim_{t \to \infty} r_0(t) = \infty$ and
such that
$(\varphi_{r_0(t)})_{t \in [1, \infty)}$ is a uniformly continuous
asymptotic homomorphism. Since $q_B \circ \varphi_{r_0(t)} = \psi_{r_0(
t)}$, we can assume that $\varphi$ is uniformly continuous. Set $t_
n = \sum_{k=1}^n \frac{1}{k}$. Let $\{V_i\}_{i=1}^{\infty}$ be a sequence
of isometries in $M(B)$ such that $\sum_{i=1}^{\infty} V_iV_i^* =
1$, in the strict topology. Set
$$
\widetilde{\nu}_t(a) = \sum_{i=2}^{\infty} V_i\varphi_{t_i + t}(a)V_i^* .
$$
It is clear that $\widetilde{\nu}: A \to M(B)$ is an asymptotic homomorphism
since $\varphi$ is. We claim that $q_B \circ \widetilde{\nu}_t$ is a completely
positive contraction for all $t \in [1,\infty)$. We prove first that $\|1-q_
B \circ \widetilde{\nu}_t(a)\| \leq 1$ when $0 \leq a \leq 1$ in $A$. To this
end, let $\{b_n\}_{n =1}^{\infty}$ be an approximate unit in $B$. Then
\begin{equation}\label{EF2}
\begin{split}
&\|1 - q_B \circ \widetilde{\nu}_t(a)\| 
\leq  \|\Bigl (1 -\sum_{k=1}^n V_kb_mV_k^*\Bigr) \Bigl( 1  - \sum_{i
=2}^{\infty} V_i\varphi_{t_i + t}(a)V_i^*\Bigr)\Bigl( 1 -\sum_{k=1}^n
V_kb_mV_k^*\Bigr)\| \\
& = \|\sum_{j > n} V_j (1 - \varphi_{t_j + t}(a))V_j^* + V_1 (1 - b_m)V_1
^* + \sum_{j =2}^n V_j(1- b_m)\left( 1  - \varphi_{t_j + t}(a) \right)(1-
 b_m)V_j^* \|
\end{split}
\end{equation}
for all $n,m$.
For a given $\epsilon > 0$ there is a $K \in \mathbb N$ such that $\|1 -
\varphi_{t_j + s}(a))\| \leq 1 + \epsilon$ for all $j \geq K$ and all $s \in
[1, \infty)$. This is because $\varphi$ is an asymptotic homomorphism. In
particular,
$$
 \|\sum_{j > K} V_j (1 - \varphi_{t_j + t}(a))V_j^*\| \leq 1 + \epsilon.
$$
We see therefore from (\ref{EF2}) that
\begin{equation}\label{EF3}
\|1 - q_B \circ \widetilde{\nu}_t(a)\| \leq \max \{ 1 + \epsilon, \  \|\sum_{j
=2}^K V_j(1- b_m)\left( 1  - \varphi_{t_j + t}(a) \right)(1- b_m)V_j^*\|
\}
\end{equation}
for all $m \in \mathbb N$. Since $q_B \circ \varphi_{t_j + t}$ is a
completely positive contraction for all $j = 2,3, \dots, K$, there is an $m \in
\mathbb N$ such that
$$
\|(1- b_m)\left( 1  - \varphi_{t_j + t}(a) \right)(1- b_m)\| \leq 1 +
\epsilon
$$
for all $j \in \{2,3,\dots, K\}$. Since $\epsilon > 0$ is arbitrary, it
follows therefore from (\ref{EF3}) that $\|1 - q_B \circ \widetilde{\nu}_t(a)\|
\leq 1$. Similar arguments show that $q_B \circ \widetilde{\nu}_t$ is a linear
self-adjoint contraction, and combined with what we have just established this
implies that $q_B \circ \widetilde{\nu}_t$ a positive linear contraction. By
applying the same argument to the maps $M_n(A) \to M_n(Q(B))$ induced by $q_B
\circ \widetilde{\nu}_t$, we see that $q_B \circ \widetilde{\nu}_t$ is a completely
positive contraction for all $t$. Let $r : [1, \infty) \to [1,\infty)$ be a
continuous function such that
\begin{equation}\label{EF7}
t_n + t \leq r(t) \leq t_{n+1} + t
\end{equation}
for $t \in [n,n+1]$. Set
\begin{equation*}
\begin{split}
S^n_t = 1 - V_nV_n^* - V_{n+1}V_{n+1}^*& +  \cos \left(\frac{\pi}{2}(t-n)
\right) V_nV_n^* +  \cos \left(\frac{\pi}{2}(t-n)\right) V_{n+1}V_{n+1}^* \\
&+  \sin \left(\frac{\pi}{2}(t-n)\right) V_{n+1}V_n^* -  \sin
\left(\frac{\pi}{2}(t-n)\right) V_{n}V_{n+1}^*
\end{split}
\end{equation*}
for $t \in [n,n+1]$. Then $S_t = S^n_tS^{n-1}_n\dots S^1_2$, $t \in [n, n+1
]$, is a norm-continuous path of unitaries in $M(B)$, constructed such that
(\ref{EF7}) and the uniform continuity of $\varphi$ ensure that
$$
\lim_{t \to \infty} S_t \Bigl( V_1 \varphi_{r(t)}(a)V_1^* + \sum_{i=2}^
{\infty} V_i \varphi_{t_{i} + t}(a)V_i^*\Bigr)S_t^* - \sum_{i = 1}^{\infty}
V_i \varphi_{t_{i+1} + t}(a)V_i^* = 0
$$
for all $a \in A$. Since $V = \sum_{i=2}^{\infty} V_iV_{i-1}^*$ is a an
isometry in $M(B)$ such that $V \left(\sum_{i = 1}^{\infty} V_i \varphi_
{t_{i+1} + t}(a)V_i^*\right)V^* = \widetilde{\nu}_t(a)$ for all $a,t$, we can
put $\mu = \nu \oplus 0$ and $U_t=S_t\oplus 1$.
\end{proof}

\begin{lem}\label{basic} Let $\psi : A \to Q(B)$ be an asymptotically split
extension. It follows that there is a asymptotically split extension $\nu : A \to Q(B)$ and a unitary $U \in M_2(M(B))
$ such that
$$
\Ad q_{M_2(B)}(U) \circ \left( \begin{smallmatrix} \psi & 0 \\ 0 & \nu
\end{smallmatrix} \right) =
\left( \begin{smallmatrix} 0 & 0 \\ 0 & \nu \end{smallmatrix} \right).
$$
\end{lem}
\begin{proof}  By assumption there is an asymptotic homomorphism
$\widetilde{\varphi} =(\widetilde{\varphi}_t)_{t \in [1,\infty)}: A \to M(B)$
such that $q_B \circ \widetilde{\varphi}_t = \psi$ for all $t \in [1,\infty)$. Set
$$
\mathcal B = \{ f \in C_b([1,\infty), M(B)) : q_B(f(t)) = q_B(f(1)), t
\in [1,\infty) \}.
$$
Then $\widetilde{\varphi}$ defines a $*$-homomorphism $\widehat{\varphi} : A
\to \mathcal B/C_0([1,\infty), B)$. By the Bartle-Graves selection theorem
there is a continuous section $\chi :  \mathcal B/C_0([1,\infty), B) \to
\mathcal B$ for the quotient map $\mathcal B \to \mathcal B/C_0([1,\infty),
B)$.
Set $\varphi_t(a) = \chi\left(\widehat{\varphi}(a)\right) (t)$. Then $\left(\varphi_t\right)_{t \in [1,\infty)}$ is an asymptotic homomorphism such that $q_B \circ \varphi_t = \psi$ for all $t$. Compared to $\widetilde{\varphi}$, $\varphi$ has
the property of being equi-continuous. This will be helpful in the
construction of $\nu$. Let $\{d_i\}_{i \in \mathbb N}$
be a dense sequence in $A$.
Let $\{\varphi_{t_i}\}_{i \in \mathbb N}$ be a
discretization of $\varphi$, cf.
\cite{MT1}, chosen such that
$\sup_{t\in[t_{i-1},t_i]}\|\varphi_t(d_k) - \varphi_{t_{i-1}}(d_k)\|
\leq \frac{1}{i}$, $k \leq i$, for all $i$.
For each $i \in \mathbb N$, choose
$\alpha_i \geq 1$ such that
\begin{equation}\label{A2}
\|\varphi_{t_i + \frac{t}{\alpha_i}}(d_k) - \varphi_{t_i}(d_k)\| \leq \frac
{1}{i}
\end{equation}
for all $t \in [1,i]$ and all $k \leq i$.
Let $\{V_i\}_{i \in \mathbb N}$ be a
sequence of isometries in $M(B)$ such that $\sum_{i =1}^{\infty} V_iV_i^*
= 1$, in the strict topology. Set
$$
\widetilde{\nu}_t(d) = \sum_{i=2}^{\infty} V_{i} \varphi_{t_i + \frac{t}
{\alpha_i}}(d)V_{i}^* .
$$
To check that $\widetilde{\nu}$ is an asymptotic homomorphism we must check that
\begin{equation}\label{A1}
\lim_{t \to \infty} \sup_{i \in \mathbb N} \|\varphi_{t_i + \frac{t}{\alpha
_i}}(a) \varphi_{t_i + \frac{t}{\alpha_i}}(b) - \varphi_{t_i + \frac{t}
{\alpha_i}}(ab)\| = 0
\end{equation}
for any pair $a,b \in A$. Let $\epsilon > 0$. Since $\varphi$ is an
asymptotic homomorphism, there is an $N \in \mathbb N$ so large that
$$
\|\varphi_{t_i + \frac{t}{\alpha_i}}(a) \varphi_{t_i + \frac{t}{\alpha_i}}(
b) - \varphi_{t_i + \frac{t}{\alpha_i}}(ab)\| \leq \epsilon
$$
for all $t \in [1,\infty)$, when $i \geq N$. Choose $T \in [1,\infty)$ so
large that
$$
\sup_{i \leq N} \|\varphi_{t_i + \frac{t}{\alpha_i}}(a) \varphi_{t_i +
\frac{t}{\alpha_i}}(b) - \varphi_{t_i + \frac{t}{\alpha_i}}(ab)\| \leq
\epsilon
$$
when $t > T$. Then
$$
\sup_{i \in \mathbb N} \|\varphi_{t_i + \frac{t}{\alpha_i}}(a) \varphi_{t_i
+ \frac{t}{\alpha_i}}(b) - \varphi_{t_i + \frac{t}{\alpha_i}}(ab)\| \leq
\epsilon
$$
when $t > T$, proving (\ref{A1}). The other asymptotic algebraic identities
follow in the same way. Since $\varphi_t(a) - \varphi_s(a) \in B$ for all
$a,s,t$, (\ref{A2}) ensures that $\widetilde{\nu}_t(d_k) - \widetilde{\nu}_1(d_k)
\in B$ for all $k$ and $t$. Thanks to the equi-continuity of $\varphi$, $q_B
\circ \widetilde{\nu}_t$ is a continuous map for each $t$, so the density of
$\{d_i\}_{i\in\mathbb N}$ in $A$ ensures that
$\nu = q_B \circ \widetilde{\nu}_t$ is independent
of $t$, and hence defines an extension which is asymptotically split by
construction. Note that $\psi \oplus \nu$ is unitarily equivalent to $\mu
= q_B \circ \widetilde{\mu}$, where
$$
\widetilde{\mu}_t(a) = V_1 \varphi_{t_1 + \frac{t}{\alpha_1}}(a)V_1^* + \sum_{
i=2}^{\infty} V_{i} \varphi_{t_i + \frac{t}{\alpha_i}}(a)V_{i}^* .
$$
Set $S = \sum_{i=2}^{\infty}V_{i-1}V_i^*$. Then $S$ is an isometry in
$M(B)$ which by (\ref{A2}) has the property that $\Ad S \circ \widetilde{\nu}_
t(d_k) - \widetilde{\mu}_t(d_k) \in B$ for all $t,k$. It follows that $\Ad q_B(
S) \circ \nu = \mu$. Hence $\psi \oplus \nu \oplus 0$ is unitarily equivalent
to $\nu \oplus 0$. Since $\nu$ is unitarily equivalent to $0 \oplus \nu$
it follows that $\psi \oplus \nu$ is unitarily equivalent to $0 \oplus \nu$,
which is the statement of the lemma.
\end{proof}

 We will say that an asymptotic extension $\varphi : A \to Q(B)$ is
\emph{semi-invertible} when there is another asymptotic extension $\psi :
A \to Q(B)$ such that $\varphi \oplus \psi$ is split (as an asymptotic
extension).

\begin{lem}\label{gomot}
Let $\varphi,\psi:A\to Q(B)$ be two semi-invertible asymptotic extensions
that are homotopic. Then there exists a split asymptotic extension $\mu$
such that $\varphi\oplus\mu$ and $\psi\oplus\mu$ are unitarily
equivalent. In fact, there is a unitary $U \in M_2(M(B))$ such that
$$
\lim_{t \to \infty} \Ad q_{M_2(B)}\left(U \right) \circ \left(
\begin{smallmatrix} \varphi_t(a) & 0 \\ 0 & \mu_t(a) \end{smallmatrix}
\right) -  \left( \begin{smallmatrix} \psi_t(a) & 0 \\ 0 & \mu_t(a)
\end{smallmatrix} \right) = 0
$$
for all $a \in A$.

\end{lem}
\begin{proof} We may assume that $\varphi$ and $\psi$ are equi-continuous.
Choose equi-continuous asymptotic extensions $\varphi_1,\psi_1 : A \to Q(B)
$ such that $\varphi \oplus \varphi_1$ and $\psi \oplus \psi_1$ are split.
Then $\varphi \oplus \psi_1$ is homotopic to the split asymptotic extension
$\psi \oplus \psi_1$, and hence $\mu_1 = \varphi \oplus \psi_1$ is split
by Lemma \ref{crux2}. Set $\mu_2 = \psi \oplus \psi_1$, and note that
$\psi \oplus \mu_1$ is unitarily equivalent to $\varphi \oplus \mu_2$. It
follows from Lemma \ref{basic} that there are split asymptotic extensions
$\nu_i$ such that $\mu_i \oplus \nu_i$ is unitarily equivalent to $\nu_i$,
$i = 1,2$. Set $\mu = \nu_1 \oplus \nu_2$.
\end{proof}

\section{Proof of the main results}

For any $C^*$-algebra $D$ we denote by $s_D$ the stabilizing
$*$-homomorphism $s_D : D \to D \otimes \mathbb K$ given in standard
notation by $s_D(d) = d \otimes e_{11}$, and we let $\beta_B : Q(B)
\otimes \mathbb K \to Q(B \otimes \mathbb K)$ be the canonical embedding.

\begin{lem}\label{w1} There is a group homomorphism
$$F : [[SA,Q(B) \otimes \mathbb K]] \to [SA, Q(B \otimes \mathbb K) \otimes
 \mathbb K]
$$
such that $F[\varphi] = [s_{Q(B\otimes \mathbb K)} \circ (\beta_B \circ
\varphi)^f]$, when $\varphi$ is an equi-continuous asymptotic homomorphism
$\varphi : SA \to Q(B) \otimes \mathbb K$.
\end{lem}
\begin{proof} We prove first that the class of $s_{Q(B\otimes \mathbb K)}
\circ (\beta_B \circ \varphi)^f$ in $[SA,  Q(B \otimes \mathbb K) \otimes
\mathbb K]$ is independent of the choices made in the construction of $(\beta
_B \circ \varphi)^f$. The philosophy underlying the approach is due to
Higson, cf. Theorem 3.4.3 in \cite{H1}. Since $[[SA,Q(B) \otimes \mathbb K]]$
is a group there is an equi-continuous asymptotic homomorphism $\psi : SA
\to Q(B) \otimes \mathbb K$ such that
$$
\left ( \begin{smallmatrix} \varphi & 0 \\ 0 & \psi \end{smallmatrix}
\right) : SA \to Q(B) \otimes \mathbb K
$$
is homotopic to $0$. It follows that
$$
\left ( \begin{smallmatrix} \beta_B \circ \varphi & 0 \\ 0 & \beta_B \circ
\psi \end{smallmatrix} \right)
$$
is homotopic to $0$, and hence split (as an asymptotic extension)
by Lemma \ref{crux2}. Then Lemma
\ref{XX} tells us that
$$
\left ( \begin{smallmatrix} \left(\beta_B \circ \varphi\right)^f & 0 \\ 0 &
 \left(\beta_B \circ \psi\right)^f \end{smallmatrix} \right)
$$
is asymptotically split. It follows then from Lemma \ref{basic} that $[s_{Q
(B \otimes \mathbb K)} \circ (\beta_B \circ \varphi)^f] = - [s_{Q(B \otimes
\mathbb K)} \circ (\beta_B \circ \psi)^f]$ in $ [SA, Q(B \otimes \mathbb K
) \otimes \mathbb K]$. Thus $[s_{Q(B \otimes \mathbb K)} \circ (\beta_B
\circ \varphi)^f] = - [s_{Q(B \otimes \mathbb K)} \circ (\beta_B \circ \psi)
^f]$, regardless of the folding data used in the construction of $(\beta_B
\circ \varphi)^f$. Assume next that $\varphi, \varphi' : SA \to Q(B) \otimes
\mathbb K$ are homotopic as asymptotic homomorphisms. By Lemma \ref{gomot}
there is a (split) asymptotic extension $\mu : SA \to Q(B \otimes \mathbb
K)$ such that $(\beta_B \circ \varphi) \oplus \mu$ and $(\beta_B \circ
\varphi') \oplus \mu$ are unitarily equivalent. Since we are now free to
choose the foldings at will we conclude from Lemma \ref{unitary} that
$[s_{Q(B \otimes \mathbb K)} \circ \left( \beta_B \circ
\varphi\right)^f] = [s_{Q (B \otimes \mathbb K)} \circ \left( \beta_B
\circ \varphi'\right)^f]$. Hence $F : [[SA,Q(B) \otimes \mathbb K]] \to
[SA, Q(B \otimes \mathbb K) \otimes \mathbb K]$ is well-defined, and by
using the freedom of choice of foldings again it follows that $F$ is a
homomorphism.
\end{proof}

Let $ i : [[SA, Q(B) \otimes \mathbb K]]_{cp} \to  [[SA, Q(B) \otimes
\mathbb K]]$ and $I : [SA, Q(B \otimes \mathbb K) \otimes \mathbb K]
\to [[SA, Q(B \otimes \mathbb K) \otimes \mathbb K]]_{cp}$ be the obvious
(forgetful) maps.
Denote by $\widehat{s_B}:Q(B)\to Q(B\otimes\mathbb K)$ the homomorphism
induced by $s_B:B\to B\otimes\mathbb K$.

\begin{lem}\label{w2} The composition
\begin{equation*}
\begin{xymatrix}{
[[SA, Q(B) \otimes \mathbb K]]_{cp} \ar[r]^-{i} & [[SA, Q(B) \otimes \mathbb
K]]  \ar[r]^-{I \circ F} & [[SA, Q(B \otimes \mathbb K) \otimes \mathbb K
]]_{cp} }
\end{xymatrix}
\end{equation*}
agrees with $(\widehat{s_B} \otimes \id_{\mathbb K})_*$.
\end{lem}
\begin{proof} Let $\varphi', \psi' : SA \to Q(B) \otimes \mathbb K$ be
completely positive asymptotic homomorphisms such that $\psi'$ represents the
inverse of $[\varphi']$ in $[[SA, Q(B) \otimes \mathbb K]]_{cp}$. Set
$\varphi = \beta_B \circ \varphi'$ and $\psi = \beta_B \circ \psi'$. It follows
from Lemma \ref{crux2} that there is an equi-continuous asymptotic homomorphism
$\pi : SA \to M_2(M(B))$ such that
\begin{equation}\label{anne}
\left( \begin{smallmatrix} \varphi_t & 0 \\ 0 & \psi_t \end{smallmatrix}
\right) = q_{M_2(B \otimes \mathbb K)} \circ \pi_t
\end{equation}
for all $t$. Let $(\widehat{\varphi}_t),
(\widehat{\psi}_t)_{t\in[1,\infty)} : SA \to M(B)$
be equi-continuous lifts of $\varphi$ and $\psi$,
respectively. Let $\left( \left(\widehat{\varphi}_t\right)_{t \in [1,\infty)}, \{u_n\}_{n=0}
^{\infty}, \{n_i\}_{i=0}^{\infty}, \{t_i\}_{i=1}^{\infty}, \right)$
be the folding data used to define $\varphi^f$, so that
$$
\varphi^f(a) = q_{B \otimes \mathbb K}\Bigl( \sum_{j=0}^{\infty} \Delta
_j\widehat{\varphi}_{t_{j+1}}(a) \Delta_j \Bigr),
$$
where $\Delta_0 = \sqrt{u_{n_0}}$ and $\Delta_j = \sqrt{u_{n_j} - u_{n_
{j-1}}}$, $j \geq 1$.  Let $F_1 \subseteq F_2 \subseteq F_3 \subseteq \dots$
be a sequence of finite sets with dense union in $SA$. Since we can choose
the folding data at will we can apply Lemma \ref{ex} to arrange that
\begin{equation}\label{AAA}
\|\Delta_j \widehat{\psi}_t(a) - \widehat{\psi}_t(a) \Delta_j \| \leq 2^{-j}
\end{equation}
for all $a \in F_j$, $t \in [1,j+3]$, and all $j$. It follows then that
$$
\psi_t(a) = q_{B\otimes \mathbb K} \circ \widehat{\psi}_t(a) = q_{B \otimes
\mathbb K} \Bigl( \sum_{j=0}^{\infty} \Delta_j \widehat{\psi}_t(a) \Delta_j
\Bigr) ,
$$
for all $t \in [1,\infty)$, first for all $a \in \bigcup_n F_n$, and then
by continuity for all $a \in SA$. We now proceed much as in the proof of
Lemma \ref{XX}: Using (\ref{anne}) we choose a sequence $r_1 < r_2 < r_3 <
\dots $ in  $\mathbb N$ such that
\begin{equation}\label{AXA}
\|\left ( \begin{smallmatrix} 1- u_j & {} \\ {} & 1- u_j \end{smallmatrix}
\right )\left[ \left ( \begin{smallmatrix} \widehat{\varphi}_{t}(d) & {} \\ {}
&  \widehat{\psi}_{t}(d)) \end{smallmatrix} \right ) - \pi_{t}(d)\right]\| \leq
 \frac{1}{n}
\end{equation}
for all $t \in [1,n+3]$, $j \geq r_{n-2}$, $d \in F_n$, and such that
(a)-(d) and (a')-(d')
from the proof of Lemma \ref{XX} hold, with $u'_j = u_j$. We
then construct the same path $\{w_i(t)\}_{i=0}^{\infty}$,
$t \in [1,\infty)$,
of unit sequences as in the proof of Lemma \ref{XX}, but this time we set
\begin{eqnarray*}
\rho_t(a) &=& \left( \begin{matrix}  \sum_{i=1}^{\infty} \Delta_i(t)
\widehat{\varphi}_{f^i(t)}(a) \Delta_i(t)   & 0  \\  0  &  \sum_{i=1}^{\infty}
\Delta_i(t) \widehat{\psi}_{t}(a) \Delta_i(t) \end{matrix} \right) \\
&
+&  \left( \begin{matrix} \Delta_0(t) &  0 \\  0   &  \Delta_0(t)
\end{matrix} \right)\pi_t(a)\left( \begin{matrix} \Delta_0(t) &  0 \\  0   &
\Delta_0(t)  \end{matrix} \right) .
\end{eqnarray*}
As in the proof of Lemma \ref{XX} we see that $\rho$ is an asymptotic
homomorphism. Since
$$
q_{B \otimes \mathbb K} \circ \rho_t = \left( \begin{smallmatrix} \varphi^f
& 0 \\ 0 & \psi_t \end{smallmatrix} \right)
$$
for all $t$, it follows from Lemma \ref{cp} that $[s_{Q(B \otimes \mathbb K
)} \circ \varphi^f] = - [s_{Q(B \otimes \mathbb K)} \circ  \psi]$ in $[[
SA, Q(B  \otimes \mathbb K) \otimes \mathbb K]]_{cp}$. Since $ - [s_{Q(B
\otimes \mathbb K)} \circ  \psi] = [s_{Q(B \otimes \mathbb K)} \circ
\varphi]$, we conclude that $I \circ F \circ i = \left( s_{Q(B \otimes
\mathbb K)} \circ \beta_B\right)_*$. Note that $ s_{Q(B \otimes \mathbb K)}
\circ\beta_B$ and $\widehat{s_B} \otimes \id_{\mathbb K}$ both map $Q(B)
\otimes \mathbb K$ into $Q(B) \otimes \mathbb K \otimes \mathbb K$, and that
there is a unitary $U \in M(Q(B) \otimes \mathbb K \otimes \mathbb K)$ such
that $ \Ad U \circ s_{Q(B \otimes \mathbb K)} \circ \beta_B =\widehat{s_B}
\otimes \id_{\mathbb K}$. Since the unitary group of $M(Q(B) \otimes \mathbb
K \otimes \mathbb K)$ is connected in the strict topology, it follows that
$s_{Q(B \otimes \mathbb K)} \circ \beta_B$ and $\widehat{s_B} \otimes
\id_{\mathbb K}$ are homotopic. Hence $I \circ F \circ i = (\widehat{s_B}
\otimes \id_{\mathbb K})_*$.
\end{proof}

Note that since $B$ is stable there is an isometry $V \in Q(B \otimes
\mathbb K)$ such that $Q(B) \ni x \mapsto V^*\widehat{s_B}(x)V \in
Q(B\otimes\mathbb K)$ is an isomorphism, which we will denote by $\gamma$,
and such that $\Ad V \circ \gamma = \widehat{s_B}$. We can therefore
improve Lemma \ref{w2} as follows.

\begin{lem}\label{w3} The composition
\begin{equation*}
\begin{xymatrix}{
[[SA, Q(B) \otimes \mathbb K]]_{cp} \ar[d]_-{i} & [[SA, Q(B \otimes \mathbb
 K) \otimes \mathbb K]]_{cp} \ar[d]^-{\left(\gamma^{-1} \otimes \id_{\mathbb
K}\right)_* } & \\
[[SA, Q(B) \otimes \mathbb K]]  \ar[ur]_-{I \circ F} &  [[SA, Q(B) \otimes
\mathbb K]]_{cp} }
\end{xymatrix}
\end{equation*}
is the identity.
\end{lem}

In a similar way we can prove

\begin{lem}\label{w4} The composition
\begin{equation*}
\begin{xymatrix}{
[[SA, Q(B) \otimes \mathbb K]] \ar[d]_-{I \circ F} & [[SA, Q(B \otimes
\mathbb K) \otimes \mathbb K]] \ar[d]^-{(\gamma^{-1} \otimes \id_{\mathbb K})_*
} & \\
[[SA, Q(B\otimes \mathbb K) \otimes \mathbb K]]_{cp}  \ar[ur]_-{i} &  [[SA,
Q(B) \otimes \mathbb K]]_{cp} }
\end{xymatrix}
\end{equation*}
is the identity.
\end{lem}

In fact, the proof is somewhat simpler because we don't have to worry about
complete positivity of the asymptotic homomorphisms.

\begin{cor}
The map $i : [[SA,Q(B) \otimes \mathbb K]]_{cp}\to [[SA,Q(B) \otimes
\mathbb K]]$ is an isomorphism.
\end{cor}
\begin{proof}
It follows from Lemma \ref{w3} that $i$
is injective. Since the $\gamma^{-1}_*
\circ i = i \circ \gamma^{-1}_*$ we conclude from Lemma \ref{w4} that $i$
is also surjective.
\end{proof}

Before we finally relate $E$-theory and $KK$-theory we need to identify
$KK(SA,Q(B))$ with $[[SA,Q(B)\otimes\mathbb K]]_{cp}$. By \cite{H-LT} and
\cite{DL} one has $KK(SA,D)=[[SA,D\otimes\mathbb K]]_{cp}$ for any separable
$C^*$-algebra $D$. To pass from separable $C^*$-algebras to $Q(B)$ we use
the following statement.

\begin{lem}\label{sep} Let $A$ and $B$ be $C^*$-algebras, $A$ separable and $B$ $\sigma$-unital.
Then
 $$
\varinjlim_{D}[[A,D\otimes\mathbb
K]]_{cp}=[[A,B\otimes\mathbb K]]_{cp},
\quad
\varinjlim_{D}KK(A,D)=KK(A,B),
 $$
where the limits are taken over the net of separable $C^*$-subalgebras
$D$ of $B$.
\end{lem}
\begin{proof} There is an obvious map
$\varinjlim_{D}[[A,D\otimes\mathbb
K]]_{cp}\to[[A,B\otimes\mathbb K]]_{cp}$. Let $\varphi : A \to B \otimes
\mathbb K$ be a completely positive asymptotic homomorphism. Let
$\{e_{ij}\}_{i,j\in \mathbb N}$ be the standard matrix units in $\mathbb K$.
Then
$$
\{s_B^{-1}\left(\left(1_B \otimes e_{1i}\right)\varphi_t(a)\left(1_B \otimes
e_{j1}\right)\right): t\in [1, \infty), i,j \in \mathbb N, a \in A\}
$$
generates a separable $C^*$-subalgebra $D$ of $B$ such that $\varphi_t(A)
\subseteq D \otimes \mathbb K$. This shows that the map is surjective.
Handling a homotopy of completely positive asymptotic homomorphisms in
the same way shows that the map is also injective. By \cite{Cuntz},
$KK(A,B) = [qA, B \otimes \mathbb K]$, and since $qA$ is separable when
$A$ is, the same argument shows that $\varinjlim_{D}KK(A,D)=KK(A,B)$.
\end{proof}

Since $Q(B)$ is unital we can combine Lemma \ref{sep} with \cite{H-LT} and
\cite{DL} to conclude that
\begin{equation}\label{ny1}
KK(SA,Q(B)) = [[SA,Q(B)\otimes \mathbb K]]_{cp}.
\end{equation}

\begin{lem}\label{ny} Let $A$ be a separable $C^*$-algebra. Then the functor
$[[SA, - \otimes \mathbb K]]$ is half-exact on the category of all
$C^*$-algebras, and there is a natural isomorphism $[[SA, Q(B)\otimes
\mathbb K]] = [[SA,SB\otimes \mathbb K]]$, for any stable $C^*$-algebra
$B$.
\end{lem}
\begin{proof} By working with equi-continuous asymptotic homomorphisms we
can adopt the proof of Lemma \ref{sep} to conclude that
\begin{equation}\label{ny5}
[[SA, B \otimes \mathbb K]] = \varinjlim_{D}[[SA,D\otimes\mathbb
K]],
\end{equation}
for any $C^*$-algebra $B$, where we take the limit over the net of separable
$C^*$-subalgebras of $B$. Therefore the fact that $[[SA, -\otimes
\mathbb K]]$ is half-exact on the category of separable $C^*$-algebras,
cf. \cite{CH} and \cite{DL}, implies that the same is true on the category
of all $C^*$-algebras. Hence the canonical extension $0 \to B \to M(B)
\to Q(B) \to 0$ gives rise to a long exact sequence, a part of which is
\begin{equation*}
\begin{xymatrix}{
[[SA, SM(B) \otimes \mathbb K]] \ar[r] &[[SA,SQ(B) \otimes \mathbb K]] \ar[r]^-{\partial} & [[SA,B \otimes \mathbb K]] \ar[r]  & [[SA,M(B)\otimes \mathbb K]],}
\end{xymatrix}
\end{equation*}
cf. e.g. Theorem 21.4.3 of \cite{Bl}. To show that $\partial$ is an
isomorphism it suffices now to show that $[[SA,SM(B)\otimes \mathbb K]] =
[[SA,M(B) \otimes \mathbb K]] = 0$. Since $B$ is stable there is a sequence
of isometries $V_i, i = 1,2,3, \dots$, in $M(B)$ such that $\psi(m) =
\sum_{i=1}^{\infty} V_imV_i^*$ converges in the strict topology for all
$m \in M(B)$ and defines an endomorphism $\psi : M(B) \to M(B)$.
Furthermore, for any pair of isometries $W_1,W_2 \in M(B)$ for which
$W_1W_1^* + W_2W_2^* = 1$,
$$
U = W_2V_1^* + \sum_{i=2}^{\infty} W_1V_{i-1}V_i^*
$$
is a unitary in $M(B)$ with the property that $W_1\psi(m)W_1^* + W_2mW_2^*
= \Ad U \circ \psi(m)$ for all $m \in M(B)$. It follows that $\psi_* + \id_*
= \psi_*$, both as endomorphisms of the group $[[SA, SM(B) \otimes
\mathbb K]]$, and of the group $[[SA, M(B)\otimes \mathbb K]]$.
Hence $[[SA, SM(B) \otimes \mathbb K]] = [[SA, M(B) \otimes \mathbb K]]
= 0$, and we conclude that $\partial$ is an isomorphism.

Applied with $A$ replaced by $SA$ we get an isomorphism
\begin{equation}\label{sq}
[[S^2A, SQ(B) \otimes \mathbb K]] = [[S^2A, B \otimes \mathbb K]].
\end{equation}
It follows from (\ref{ny5}) and \cite{DL} that the suspension maps
\begin{equation}\label{sq1}
S : [[SA, Q(B) \otimes \mathbb K]] \to [[S^2A, SQ(B)\otimes \mathbb K]]
\end{equation}
and
\begin{equation}\label{sq2}
S : [[S^2A, B \otimes \mathbb K]] \to [[S^3A, SB\otimes \mathbb K]]
\end{equation}
are both isomorphisms. Since $S^3A \otimes \mathbb K$ and $SA \otimes
\mathbb K$ are equivalent in E-theory, there is also an isomorphism
\begin{equation}\label{sq3}
[[S^3A, SB \otimes \mathbb K]] = [[SA, SB\otimes \mathbb K]].
\end{equation}
The desired isomorphism $[[SA,Q(B) \otimes \mathbb K]] = [[SA, SB \otimes
\mathbb K]]$ is put together by (\ref{sq1}), (\ref{sq}), (\ref{sq2}) and
(\ref{sq3}).

\end{proof}

\begin{thm}\label{THM} Let $A$ be a separable $C^*$-algebra, and $B$ a
stable $\sigma$-unital $C^*$-algebra. There is then a natural isomorphism
$$
E(A,B) = KK(SA, Q(B)).
$$
\end{thm}
\begin{proof} It follows from Lemma \ref{ny} that $E(A,B) = [[SA, Q(B) \otimes
\mathbb K]]$. We have already seen that
$KK(SA,Q(B)) = [[SA, Q(B) \otimes \mathbb K]]_{cp}$, cf. (\ref{ny1}).
Combine Lemma \ref{w2} and Lemma \ref{w3}.
\end{proof}

\begin{cor}\label{cor}  Let $A$ be a separable $C^*$-algebra, and $B$ a
stable $\sigma$-unital $C^*$-algebra. Then $E(A,B) = [SA, Q(B) \otimes
\mathbb K]$.
\end{cor}
\begin{proof} It follows from Lemma \ref{w3} that $I : [SA, Q(B \otimes
\mathbb K) \otimes \mathbb K] \to  [[SA, Q(B \otimes \mathbb K) \otimes
\mathbb K]]_{cp}$ is surjective. To show that it is also injective, observe that
Lemma 4.3 of \cite{MT3} tells us that the composition
\begin{equation*}
\begin{xymatrix}{
[SA, Q(B) \otimes \mathbb K] \ar[d]_-{i \circ I} & [SA, Q(B \otimes \mathbb
 K) \otimes \mathbb K] \ar[d]^-{\left(\gamma^{-1} \otimes \id_{\mathbb K}
\right)_* } & \\
[[SA, Q(B) \otimes \mathbb K]]  \ar[ur]_-{F} &  [SA, Q(B) \otimes \mathbb K
] }
\end{xymatrix}
\end{equation*}
is the identity.
\end{proof}



\vspace{2cm}
\parbox{7cm}{V. M. Manuilov\\
Dept. of Mech. and Math.,\\
Moscow State University,\\
Moscow, 119992, Russia\\
e-mail: manuilov@mech.math.msu.su
}
\hfill
\parbox{6cm}{K. Thomsen\\
Institut for matematiske fag,\\
Ny Munkegade, 8000 Aarhus C,\\
Denmark\\
e-mail: matkt@imf.au.dk
}


\begin{thebibliography}{9999}


\bibitem[A]{A}
{W. Arveson,} {\it Notes on extensions of $C^*$-algebras}.
Duke Math. J. {\bf 44} (1977), 329-355.


\bibitem[Bl]{Bl} B. Blackadar, {\em $K$-theory for Operator Algebras},
Math. Sci. Res. Inst. Publ. {\bf 5}, Springer-Verlag, New York, 1986.


\bibitem[CH]{CH} A. Connes, N. Higson, {\em D\'eformations, morphismes
asymptotiques et
$K$-th\'eorie bivariante}, C. R. Acad. Sci. Paris S\'er. I Math. {\bf 311}
(1990), 101-106.


\bibitem[C]{Cuntz}
J. Cuntz, {\em A new look at $KK$-theory}, $K$-Theory {\bf 1} (1987),
31--51.


\bibitem[DL]{DL} M. Dadarlat and T. Loring, {\em K-homology, Asymptotic
Representations, and Unsuspended E-theory}, J. Func. Analysis {\bf 126}
(1994), 367-383.



\bibitem[H1]{H1} N. Higson, {\em Algebraic $K$-theory of stable
$C\sp *$-algebras}, Adv. in Math. {\bf 67} (1988), no. 1, 140 pp.

\bibitem[H2]{H2} \bysame, {\em Categories of fractions and excision in
$KK$-theory} J. Pure Appl. Algebra {\bf 65} (1990), no. 2, 119-138.




\bibitem[H-LT]{H-LT} T. Houghton-Larsen, K. Thomsen, {\em  Universal
(co)homology theories}, K-theory {\bf 16} (1999), 1-27.




\bibitem[K]{K} G. Kasparov, {\em Equivariant KK-theory and the Novikov
conjecture}, Invent. Math. {\bf 91} (1988), 513-572.





\bibitem[L]{L} T. Loring, {\em Almost multiplicative maps between
$C^*$-algebras}, Operator Algebras and Quantum Field Theory, Rome 1996.


\bibitem[MT1]{MT1}
{V. M. Manuilov, K. Thomsen,} {\it Quasidiagonal extensions and sequentially
trivial asymptotic homomorphisms}. Adv. Math. {\bf 154} (2000), 258--279.

\bibitem[MT3]{MT3}
{V. M. Manuilov, K. Thomsen,} {\it The Connes--Higson construction is an
isomorphism}. Preprint.

\bibitem[M]{M} {V. M. Manuilov,} {\it Asymptotic homomorphisms into the
Calkin algebra}.
J. Reine Angew. Math., to appear.


\bibitem[P]{P}
{G. K. Pedersen,} {\it $C^*$-algebras and their automorphism groups}.
Academic Press, London -- New York -- San Francisco, 1979.






\bibitem[S]{S} G. Skandalis, {\em Le bifuncteur de Kasparov n'est pas exact},
C.R. Acad. Sci. Paris, S\'er. I Math. {\bf 313} (1991), 939-941.






\bibitem[T]{T}  {K. Thomsen}, {\it Homotopy invariance in E-theory},
Preprint, Aarhus, 2001.


\bibitem[V]{V}
{D. Voiculescu,} {\it A note on quasidiagonal $C^*$-algebras and homotopy}.
Duke Math. J. {\bf 62} (1991), 267--271.

\end{thebibliography}
\end{document}